\documentclass[12pt]{amsart}
\usepackage{amssymb,latexsym}

\numberwithin{equation}{section}

\newtheorem{theorem}{Theorem}[section]
\newtheorem{proposition}[theorem]{Proposition}

\newtheorem{corollary}[theorem]{Corollary}
\newtheorem{lemma}[theorem]{Lemma}

\theoremstyle{definition}

\newtheorem{remark}[theorem]{Remark}
\newtheorem{example}[theorem]{Example}
\newtheorem{definition}[theorem]{Definition}

\newtheorem{problem}[theorem]{Problem}

\newcommand{\doublesubscript}[3]{
\displaystyle\mathop{\displaystyle #1_{#2}}_{#3}}

\newcommand{\overunder}[2]{
\!\begin{array}{c}
\scriptstyle{#1}\\[-.1in]
-\!\!\!-\!\!\!-\\[-.1in]
\scriptstyle{#2}
\end{array}
\!
}

\renewcommand{\eqref}[1]{{\rm (\ref{#1})}}

\def\endproof{\hfill$\square$\medskip}

\def\AA{\mathcal{A}}

\def\ZZ{\mathbb{Z}}
\def\CC{\mathbb{C}}
\def\PP{\mathbb{P}}
\def\RR{\mathbb{R}}
\def\QQ{\mathbb{Q}}
\def\TT{\mathbb{T}}
\def\xx{\mathbf{x}}
\def\pp{\mathbf{p}}
\def\FFcal{\mathcal{F}}
\def\ZP{\ZZ\PP}

\def\Upper{\mathcal{U}}
\def\Lower{\mathcal{L}}
\def\Supp{\operatorname{Supp}}

\newcommand{\mat}[4]{\left(\!\!\begin{array}{cc}
#1 & #2 \\ #3 & #4 \\
\end{array}\!\!\right)}

\def\veps{\varepsilon}
\def\veps{\varepsilon}
\def\gg{\mathfrak{g}}
\def\hh{\mathfrak{h}}
\def\ii{\mathbf{i}}
\def\ee{\mathbf{e}}
\def\sgn{\operatorname{sgn}}

\begin{document}

\title[Cluster algebras III]
{Cluster algebras III:\\ Upper bounds and double Bruhat cells}

\author{Arkady Berenstein}
\address{Department of Mathematics, University of
Oregon, Eugene, OR 97403, USA} \email{arkadiy@math.uoregon.edu}

\author{Sergey Fomin}
\address{Department of Mathematics, University of Michigan,
Ann Arbor, MI 48109, USA} \email{fomin@umich.edu}

\author{Andrei Zelevinsky}
\address{\noindent Department of Mathematics, Northeastern University,
 Boston, MA 02115, USA}
\email{andrei@neu.edu}

\begin{abstract}
We develop a new approach to cluster algebras based on the notion of
an upper cluster algebra, defined as an intersection of 
Laurent polynomial rings. Strengthening the Laurent phenomenon
established in~\cite{fz-clust1},
we show that, under an assumption of ``acyclicity", a
cluster algebra coincides with its ``upper" counterpart,
and is finitely generated; in this case, we also describe its defining
ideal, and construct a standard monomial basis.
We prove that the coordinate ring of any double Bruhat cell in a
semisimple complex Lie group is naturally isomorphic to an
upper cluster algebra explicitly defined in terms of relevant
combinatorial data.
\end{abstract}

\date{July 4, 2003. Revised January~7, 2004}

 \thanks{Research 
supported 
by NSF (DMS) grants \# 0102382 (A.B.), 0070685 (S.F.), and 0200299~(A.Z.).}


\maketitle


\tableofcontents


\section*{Introduction}

The study of cluster algebras began 
in~\cite{fz-clust1} and continued in~\cite{fz-clust2}.
The present paper, the third in the series,
can be read independently of its predecessors.

Cluster algebras are a class of
commutative rings defined axiomatically in terms of a
distinguished family of generators called cluster variables.
The main motivation for cluster algebras came from the study of dual canonical
bases (see~\cite{z-nato} and references therein)
and from the theory of double Bruhat cells in complex semisimple Lie
groups and related total positivity criteria (see \cite{fz-double, ssvz, z-imrn}).

In this paper, we develop a new approach to cluster algebras.
The cornerstone of this approach is the new notion of an
\emph{upper cluster algebra}, defined as an intersection of certain
(infinitely many) rings of Laurent polynomials.
This notion is motivated in part by the ``Laurent phenomenon'' established
in~\cite{fz-clust1}, which can be rephrased as saying that any cluster
algebra is contained in its ``upper'' counterpart.


In Corollary~\ref{cor:upper-bound-universal}, we show
that, under mild technical assumptions, an upper cluster algebra
can be represented as a
specific \emph{finite} intersection of Laurent polynomial rings.
Our main structural result 
is Theorem~\ref{th:acyclic-sufficient-upper}:
under an additional assumption of ``acyclicity,''
we prove that a cluster algebra~$\AA$ coincides with the corresponding
upper cluster algebra.
We also show that under these assumptions,
$\AA$~is finitely generated;
we describe its defining ideal in concrete terms
and construct a standard monomial basis.
In addition, we completely classify finitely generated
skew-sym\-metriz\-able cluster algebras of rank~3.

The double Bruhat cells $G^{u,v}$ are intersections of cells taken
from the Bruhat decompositions
with respect to two opposite Borel subgroups; thus, they are associated
with pairs of elements $(u,v)$ of the Weyl group.
Double Bruhat cells make a natural appearance in the context of
quantum groups~\cite{CC}, total positivity~\cite{lusztig-tp},
and integrable systems~\cite{reshet,kogzel}.
The connection between cluster algebras and double Bruhat
cells was hinted upon in the papers~\cite{fz-clust1,fz-clust2}, which
were devoted to the study of structural properties of cluster algebras.
Establishing this connection is one of the main purposes of this paper.
In Theorem~\ref{double-cell=upper-bound-explicit},
we prove that the coordinate ring of any double Bruhat cell~$G^{u,v}$ is
naturally isomorphic to the complexification of an upper cluster algebra
explicitly defined in terms of combinatorial data associated with
$u$ and~$v$.

The paper is organized as follows.
Section~\ref{sec:main-results} presents our main results concerning the
general framework introduced in this paper.
In Section~\ref{sec:double-cells}, these results and constructions are
applied to 
double Bruhat cells and related (upper) cluster algebras.
The rest of the paper is devoted to the proofs of the statements in
Section~\ref{sec:main-results}.


\section{Upper and lower bounds}
\label{sec:main-results}

\subsection{The upper bounds}

We set up the stage by introducing some terminology and notation.
As in \cite{fz-clust1}, let $\PP$ be the \emph{coefficient group},
an abelian group without torsion, written multiplicatively.
A prototypical example is a free abelian group of finite rank.
Let $n$ be a positive integer, and let $[1,n]$ stand for the set
$\{1, \dots, n\}$.
As an \emph{ambient field} for our algebra, we take a field $\FFcal$
isomorphic to the field of rational functions in $n$ independent
variables
with coefficients in (the field of fractions of) the integer group ring~$\ZZ \PP$.
Slightly modifying the definition\footnote{
For the record, we articulate the differences with~\cite{fz-clust2}.
First, here we do not need the ``auxiliary addition''~$\oplus$
making~$\PP$ into a semifield. We also drop the
``normalization'' requirement $p_i^+\oplus p_i^-=1$,
which is the only place in the definition of a seed where the
auxiliary addition is used.
Finally, there is a notational difference:
in~\cite{fz-clust2}, we avoid explicit labeling of the
elements of $\xx$ and~$\pp$, and of the rows and columns of~$B$, by
integers from~$[1,n]$.}
in \cite[Section~1.2]{fz-clust2}, we define
a \emph{seed} in $\FFcal$ as a triple $\Sigma = (\xx, \pp, B)$, where
\begin{itemize}

\item $\xx\! =\! \{x_1, \dots, x_n\} 
$ is a
transcendence basis of $\FFcal$ over the field of fractions of~$\ZZ \PP$.

\item $\pp = (p_1^\pm, \dots, p_n^\pm)$ is a $2n$-tuple of elements of $\PP$.

\item $B \!=\! (b_{ij})$ is a \emph{sign-skew-symmetric} $n\!\times\! n$ integer
matrix \cite[Definition~4.1]{fz-clust1}, i.e.,
\begin{equation}
\label{eq:sss}
\text{
for any $i$ and $j$, either $b_{ij} = b_{ji} = 0$, or
$b_{ij} b_{ji} < 0$.
}
\end{equation}
\end{itemize}

We refer to $\xx$ as the \emph{cluster},
to its elements $x_i$ as \emph{cluster variables},
to $\pp$ as the \emph{coefficient tuple}, and to $B$ as the
\emph{exchange matrix} of a seed~$\Sigma$.
Seeds obtained from each other by a simultaneous
relabeling of the cluster variables~$x_j$,
the coefficients~$p_j^\pm$, and the matrix entries $b_{ij}$
are regarded as identical.

For $j \in [1,n]$, we define the adjacent cluster 
$\xx_j$ by
\begin{equation}
\label{eq:cluster-adjacent}
\xx_j = \xx - \{x_j\} \cup \{x'_j\},
\end{equation}
where the cluster variables $x_j$ and $x'_j$ are related by the exchange relation
\begin{equation}
\label{eq:exchange}
x_j x'_j = P_j(\xx) = p^+_j \prod\limits_{b_{ij}>0}
x_i^{b_{ij}}+ p^-_j \prod\limits_{b_{ij}<0} x_i^{-b_{ij}}.
\end{equation}
Clearly, each $\xx_j$ is also a transcendence basis of $\FFcal$.
Let $\ZP [\xx] \!=\! \ZP [x_1, \dots, x_n]$
(resp.,~$\ZP [\xx^{\pm 1}] \!=\! \ZP [x_1^{\pm 1}, \dots, x_n^{\pm 1}]$)
denote the ring~of polynomials (resp.,~Laurent polynomials) in
$x_1, \dots, x_n$ with coefficients in~$\ZP$.

\begin{definition}
\label{def:upper-bound}
For a seed~$\Sigma$, we denote by $\Upper (\Sigma)$
the $\ZP$-subalgebra of $\FFcal$ given by
\[
\Upper (\Sigma) = \ZP [\xx^{\pm 1}] \cap
 \ZP [\xx_1^{\pm 1}] \cap \cdots \cap \ZP [\xx_n^{\pm 1}]\, .
\]
Thus, $\Upper (\Sigma)$ is formed by the elements of $\FFcal$
which are expressed as Laurent polynomials over $\ZP$ in the variables from
each individual cluster $\xx, \xx_1, \cdots, \xx_n$.
Informally, we call $\Upper (\Sigma)$ the 
\emph{upper bound} associated with the seed~$\Sigma$,
since it contains the corresponding cluster algebra
(see Definition~\ref{def:ca-thru-lower-bound})
in case the latter is well-defined.
\end{definition}

Definition~\ref{def:upper-bound} makes it easy to check
whether a given element of the ambient field $\FFcal$ belongs to
$\Upper (\Sigma)$; however, it seems to be a 
difficult problem to construct sufficiently many elements
of $\Upper (\Sigma)$ with desirable properties\footnote{
Here, an optimistic interpretation of ``sufficiently many" is a
$\ZZ \PP$-linear basis of $\Upper (\Sigma)$,
while an important ``desirable property,"
in the case of double Bruhat cells, is the property
of a regular function to evaluate positively at all totally
nonnegative elements (in the sense of Lusztig;
see \cite{lusztig-tp} or \cite[Section~1.3]{fz-double}).
The latter is a crucial property that the
elements of ``dual canonical bases'' must possess, according to
Lusztig~\cite{lusztig-tp}.}.
The cluster algebra machinery 
introduced in \cite{fz-clust1,fz-clust2}
provides such a construction in terms of a recursive
process involving \emph{seed mutations}.

\subsection{Seed mutations and invariance of upper bounds}
Definitions~\ref{def:matrix-mutation} and \ref{def:seed-mutation}
below are slight modifications of \cite[Definition~4.2]{fz-clust1} and
\cite[Definition~1.1]{fz-clust2}, \hbox{respectively.}

\begin{definition}
\label{def:matrix-mutation}
Let 
$B \!=\! (b_{ij})
$
and $B' \!=\! (b'_{ij})
$
be 
real square matrices of the same size.
We say that $B'$ is obtained from $B$ by a \emph{matrix mutation} in
direction~$k$ (denoted $B' = \mu_k (B)$)
if
\begin{equation}
\label{eq:mutation}
b'_{ij} =
\begin{cases}
-b_{ij} & \text{if $i=k$ or $j=k$;} \\[.05in]
b_{ij} + \displaystyle\frac{|b_{ik}| b_{kj} +
b_{ik} |b_{kj}|}{2} & \text{otherwise.}
\end{cases}
\end{equation}
\end{definition}

It is easy to check that $\mu_k(\mu_k(B))=B$.

\begin{definition}
\label{def:seed-mutation}
Let $\Sigma = (\xx, \pp, B)$ be a seed in $\FFcal$, as above,
and let $k \in [1,n]$.
We~say that a seed
$\Sigma' = (\xx', \pp', B')$
is obtained from $\Sigma$ by a \emph{seed mutation} in direction~$k$ if
\begin{itemize}

\item
$\xx' = \xx_k$ is given by
(\ref{eq:cluster-adjacent})--(\ref{eq:exchange}).
\item
the $2n$-tuple
$\pp'=({p'}^ \pm_1, \dots, {p'}^ \pm_n)$
satisfies ${p'}^ \pm_k = p^\mp_k$ and, for $i \neq k$,
 \begin{equation}
\label{eq:p-mutation}
{{p'} ^+_i}/{{p'} ^-_i} =
\begin{cases}
(p_k^+)^{b_{ki}}{p ^+_i}/{p ^-_i}  & \text{if $b_{ki}\geq 0$};\\[.05in]
(p_k^-)^{b_{ki}}{p ^+_i}/{p ^-_i}  & \text{if $b_{ki}\leq 0$}.
\end{cases}
\end{equation}

\item
$B'$ is obtained from $B$ by the matrix
mutation in direction~$k$ (as in~\eqref{eq:mutation}).

\end{itemize}
\end{definition}

We note that a seed mutation in direction $k$ might not exist:
we require $\Sigma'$ to be a seed, which in turn demands that
the matrix given by (\ref{eq:mutation}) is sign-skew-symmetric.

In contrast with~\cite[Definition~1.1]{fz-clust2}
(but in implicit agreement with~\cite{fz-clust1}),
the seed $\Sigma'$ in Definition~\ref{def:seed-mutation}
is not unique\footnote{The uniqueness in \cite[Definition~1.1]{fz-clust2} was
achieved by supplying $\PP$ with a semifield structure and
imposing the normalization condition, but we do not need it here.},
since there are $n-1$ degrees
of freedom in the choice of a coefficient tuple $\pp'$
satisfying~(\ref{eq:p-mutation}).
This non-uniqueness is however inconsequential for our constructions:
it is straightforward to verify that, for any fixed sequence of indices
$k_1,\dots,k_N$,
the upper bound $\Upper(\Sigma')$---as well as its
``lower bound'' counterpart to be defined below---does not depend
on a particular choice of a seed $\Sigma'$ obtained from $\Sigma$
by a sequence of mutations in directions $k_1,\dots,k_N$.

If $\Sigma'$ is obtained from $\Sigma$ by a seed mutation in
direction~$k$, then, conversely, $\Sigma$ can be obtained from
$\Sigma'$ by a seed mutation in
the same direction.

Our first main result provides a justification for
Definition~\ref{def:seed-mutation} by showing that, under
mild restrictions, the algebra $\Upper(\Sigma)$ is invariant
under seed mutations.

\begin{definition}
\label{def:coprime-seed}
A seed $\Sigma$ is called \emph{coprime} if
the polynomials $P_1, \dots, P_n$ appearing in the exchange
relations (\ref{eq:exchange}) are pairwise coprime
in $\ZP [\xx]$.
That is, any common divisor of $P_i$ and $P_j$ must belong to~$\PP$.
\end{definition}

\begin{theorem}
\label{th:upper-bound-universal}
Assume that two seeds $\Sigma$ and $\Sigma'$ are related by a seed
mutation, and are both 
coprime.
Then the corresponding upper bounds coincide:
$\Upper(\Sigma) = \Upper(\Sigma')$.
\end{theorem}

Matrix mutations (resp., seed mutations) give rise to an
equivalence relation among square matrices
(resp., seeds), which we call \emph{mutation equivalence} and denote
by $B' \sim B$ (resp., $\Sigma' \sim \Sigma$).

We say that a seed $\Sigma= (\xx, \pp, B)$ is \emph{totally mutable} if it admits
unlimited mutations in all directions, i.e., if every
matrix mutation equivalent to $B$ is
sign-skew-symmetric; in this case,
$B$ is called \emph{totally sign-skew-symmetric}.
Although the latter condition is hard to verify in a concrete setting,
there is fortunately a simple sufficient condition that holds in all of
our geometric applications.
Following \cite[Definition~4.4]{fz-clust1}, we call a square matrix $B$
\emph{skew-symmetrizable} if there exists a
diagonal 
matrix $D$ with positive 
diagonal entries $d_i$ such that $DB$ is
skew-symmetric, i.e., $d_i b_{ij} = - d_j b_{ji}$ for all $i$ and~$j$.
It is easy to show \cite[Proposition~4.5]{fz-clust1}
that any skew-symmetrizable matrix is totally
sign-skew-symmetric, so any seed containing it is totally mutable.

\begin{definition}
Let $\,\Sigma_\circ\,$ be a totally mutable seed.
The \emph{upper cluster algebra}
\hbox{$\overline\AA=\overline\AA(\Sigma_\circ)$}
defined by~$\Sigma_\circ$
is the intersection of the subalgebras $\Upper (\Sigma)\subset\FFcal$ for all seeds
$\Sigma \sim \Sigma_\circ$.
In other words, $\overline\AA$
consists of the elements of $\FFcal$
which are Laurent polynomials over~$\ZP$
in the cluster variables from any given seed that is mutation equivalent
to~$\Sigma_\circ$.
\end{definition}

Theorem~\ref{th:upper-bound-universal} has the following direct implication.

\begin{corollary}
\label{cor:upper-bound-universal}
If all seeds mutation equivalent to a totally mutable
seed $\Sigma_\circ$ are coprime,
then the upper bound $\Upper(\Sigma)$ is independent of the choice
of~$\Sigma\sim\Sigma_\circ$, and so is equal to the upper cluster
algebra~$\overline\AA(\Sigma_\circ)$.
\end{corollary}


The coprimality condition in
Corollary~\ref{cor:upper-bound-universal} is satisfied for a
large and important class of cluster algebras of \emph{geometric type}
\cite[Definition~5.7]{fz-clust1}.
For these algebras, the coefficient group $\PP$ is a multiplicative free
abelian group; we denote its generators by $x_{n+1},
\dots, x_m$, for some $m \geq n$.
For a seed $\Sigma$ of geometric type, the coefficient tuple
is encoded by an
$(m-n) \times n$ integer matrix $B^c = (b_{i,j})_{\substack{n < i \leq
  m\\ 1 \leq j \leq n}}$
as follows:
$$
p_j^+ = \prod_{\substack{n < i \leq m\\ b_{ij}>0}} x_i^{b_{ij}},
\qquad
p_j^- = \prod_{\substack{n < i \leq m\\ b_{ij}<0}}
x_i^{-b_{ij}} .
$$
Thus, a totally mutable seed $\Sigma$ of geometric type can be thought
of as a pair $(\xx, \tilde B)$
consisting of a cluster $\xx$ and an $m \times n$ matrix $\tilde B$
with blocks $B$ and $B^c$ such that the \emph{principal part} $B$
of $\tilde B$ is totally sign-skew-symmetric.
Furthermore, the seed mutation in each direction $k$ is normalized so that the matrix
$\tilde B'$ of the mutated seed $\Sigma'$ is obtained from $\tilde B$
by extending the matrix mutation rule (\ref{eq:mutation}) to all
$i \in [1,m]$ and $j \in [1,n]$.
(One checks that this rule is compatible with
(\ref{eq:p-mutation}).)

\begin{proposition}
\label{pr:coprime}
Let $\Sigma = (\xx, \tilde B)$ be a seed of 
geometric type. 
%
If the matrix $\tilde B$ has full rank, then all seeds mutation equivalent to
$\Sigma$ are coprime.
\end{proposition}

Combining Corollary~\ref{cor:upper-bound-universal} with
Proposition~\ref{pr:coprime}, we obtain the following result.

\begin{corollary}
\label{cor:upper-bound-from-matrix}
Let $\tilde B$ be an $m \times n$ 
integer matrix of rank~$n$
whose principal part is skew-symmetrizable,
and let $\Sigma = (\xx, \tilde B)$ be the corresponding (totally
mutable) seed of geometric type.
Then the upper bound $\Upper(\Sigma)$ depends solely on the mutation
equivalence class of~$\Sigma$,
and so coincides with the upper cluster algebra
$\overline\AA=\overline\AA(\Sigma)$.
\end{corollary}


In Section~\ref{sec:double-cells},
Corollary~\ref{cor:upper-bound-from-matrix} is applied to the study of
coordinate rings of double Bruhat cells
$G^{u,v}$ in a complex semisimple Lie group~$G$.
We show that any such ring is isomorphic
to an upper cluster algebra~$\overline\AA$ which is
constructed as in Corollary~\ref{cor:upper-bound-from-matrix} from
a matrix~$\tilde B$ determined from the Weyl group elements $u$
and~$v$ by an explicit combinatorial procedure.

\pagebreak[3]

\vbox{
\subsection{Lower bounds and cluster algebras}

\begin{definition}
\label{def:lower-bound}
The \emph{lower bound} $\Lower (\Sigma)$
associated with a seed $\Sigma$ is
defined by
\[
\Lower (\Sigma) = \ZP[x_1, x'_1, \dots, x_n, x'_n] \,,
\]
where $x_1,\dots,x_n$ and $x'_1,\dots,x'_n$ have the same meaning as
in \eqref{eq:cluster-adjacent}--\eqref{eq:exchange}.
Thus, $\Lower (\Sigma)$ is the $\ZP$-subalgebra of $\FFcal$ generated by the union
of $n\!+\!1$ clusters $\xx, \xx_1, \dots, \xx_n$.
\end{definition}
}

\begin{definition}
\label{def:ca-thru-lower-bound}
The \emph{cluster algebra} $\AA = \AA(\Sigma_\circ)$
associated with a totally mutable seed~$\Sigma_\circ$ is the $\ZP$-subalgebra of
$\FFcal$ generated by the union of all lower bounds $\Lower (\Sigma)$
for $\Sigma \sim \Sigma_\circ$.
\end{definition}

In other words, $\AA$ is the $\ZP$-subalgebra of
$\FFcal$ generated by all cluster variables from all seeds $\Sigma$
mutation equivalent to~$\Sigma_\circ$.
This definition is consistent with the ones in~\cite{fz-clust1}
and~\cite{fz-clust2}, with one important distinction: the ground ring for
our cluster algebra~$\AA$ is~$\ZP$, the integral group
ring of the coefficient group~$\PP$.
(In~\cite{fz-clust2}, the ground ring was the subring of~$\ZP$
generated by the coefficients of all exchange relations,
while in~\cite{fz-clust1}, we allowed any subring of~$\ZP$ containing
these coefficients.)

The above results imply the following corollary that
justifies our use of the
``lower bound''/``upper bound'' terminology.

\begin{corollary}
\label{cor:three-inclusions}
Let $\Sigma_\circ$ be a totally mutable seed.
Then, for any $\Sigma\sim\Sigma_\circ$, we have
\begin{equation}
\label{eq:three-inclusions}
\Lower(\Sigma)\subset\AA(\Sigma_\circ)\subset\overline\AA(\Sigma_\circ)\subset\Upper(\Sigma).
\end{equation}
\end{corollary}

\begin{proof}
First suppose that any $\Sigma\sim \Sigma_\circ$ is coprime.
Then Corollary~\ref{cor:upper-bound-universal} applies,
and we obtain
$\overline\AA(\Sigma_\circ)=\Upper(\Sigma)$.
On the other hand, it is obvious from the exchange relation
\eqref{eq:exchange}
that $\Lower (\Sigma) \subset  \Upper(\Sigma)$, implying
\eqref{eq:three-inclusions} (with the last inclusion becoming an
equality).
The coprimality assumption can be lifted by using the same argument as
in~\cite[page~507]{fz-clust1}.
\end{proof}

\begin{remark}
\label{rem:new-proof-Laurent}
The inclusion $\AA(\Sigma_\circ)\subset\overline\AA(\Sigma_\circ)$ in
\eqref{eq:three-inclusions} codifies the
\emph{Laurent phenomenon} first proved in
\cite[Theorem~3.1]{fz-clust1}: every cluster variable is expressed in
terms of any
given cluster as a Laurent polynomial with coefficients in~$\ZP$.
Our proof of Theorem~\ref{th:upper-bound-universal} given in
Section~\ref{sec:proof-upper-bound} does not rely on
\cite[Theorem~3.1]{fz-clust1},
so we obtain a new proof of the Laurent phenomenon. 
\end{remark}


\subsection{Acyclicity and closing of the gaps}

The structure of the upper bounds $\Upper (\Sigma)$
(or an upper cluster algebra~$\overline\AA$)
is not well understood in general.
On the positive side, we will now formulate an extra condition on~$\Sigma$,
which holds in many important applications,
and makes the structure of $\Upper (\Sigma)$ much more tractable.

Recall that the exchange matrix $B$ of a seed $\Sigma$ is
sign-skew-symmetric, i.e., it satisfies~(\ref{eq:sss}).
We encode the sign pattern of matrix entries of $B$ by the
directed graph $\Gamma (\Sigma) = \Gamma (B)$ with the vertices $1, \dots, n$ and
the directed edges $(i,j)$ for $b_{ij} > 0$.

\begin{definition}
\label{def:acyclic-seed}
A seed $\Sigma$ (and the corresponding exchange matrix~$B$) is called
\emph{acyclic} if 
$\Gamma (\Sigma)$ has
no oriented cycles. 
\end{definition}


We first investigate the relevance of the acyclicity condition for
the structure of the lower bound~$\Lower(\Sigma)$.

\begin{definition}
\label{def:standard-monomials}
Let $\Sigma$ be a seed, and let the elements $x_1,\dots,x_n$ and
$x'_1,\dots,x'_n$ have the same meaning as
in \eqref{eq:cluster-adjacent}--\eqref{eq:exchange}.
A~monomial in $x_1, x'_1, \ldots, x_n, x'_n$ is called
\emph{standard} if it contains no product of the form~$x_j x'_j$.
\end{definition}

It is routine to show that the standard monomials
always span $\Lower(\Sigma)$ as a $\ZZ \PP$-module:
the relations (\ref{eq:exchange}) can be viewed as straightening
relations for a term order in which $x_1',\dots,x_n'$ are much more expensive than
$x_1,\dots,x_n\,$, so that the
monomial on the left-hand side of (\ref{eq:exchange})
is more expensive than all the monomials
on the right-hand side.


\begin{theorem}
\label{th:acyclic-sufficient-lower}
Let $\Lower(\Sigma)=\ZP[x_1,x'_1,\dots,x_n,x'_n]$
be the lower bound associated with a seed~$\Sigma$.
The standard monomials in $x_1, x'_1, \ldots,
x_n, x'_n$ are linearly independent over~$\ZP$
(that is, they form a $\ZP$-basis of $\Lower(\Sigma)$)
if and only if $\Sigma$ is acyclic.
\end{theorem}

Theorem~\ref{th:acyclic-sufficient-lower} has the following
immediate implication.

\begin{corollary}
\label{cor:Grobner}
Let $\Lower (\Sigma)$
be the lower cluster algebra associated with a seed $\Sigma$,
and let $\mathcal{I}$ denote the ideal of relations
among its generators $x_1,\dots,x_n,x_1',\dots,x_n'$. \linebreak[2]
If $\Sigma$ is acyclic, then
the polynomials $x_j x'_j - P_j(\xx)$, for $j \in [1,n]$, generate the
ideal~$\mathcal{I}$.
Moreover, these polynomials form a Gr\"obner basis
for~$\mathcal{I}$ with respect to any term order
in which $x_1',\dots,x_n'$ are much more expensive than $x_1,\dots,x_n\,$.
\end{corollary}

Our next result shows that the acyclicity condition closes the gap
between the upper and lower bounds.

\begin{theorem}
\label{th:acyclic-sufficient-upper}
If a seed $\Sigma$ is coprime and acyclic, then
$\Lower (\Sigma) 
= \Upper (\Sigma)$.
\end{theorem}


\begin{corollary}
\label{cor:finite-generation-sufficient}
If a cluster algebra possesses a coprime and acyclic seed,
then it coincides with
the corresponding upper cluster algebra.
\end{corollary}

Our next result is based on
Theorem~\ref{th:acyclic-sufficient-upper};
however, it does not require coprimality.

\begin{theorem}
\label{th:lower-bound-exact}
The cluster algebra $\AA(\Sigma)$
associated with a totally mutable seed $\Sigma$
is equal to the lower bound $\Lower(\Sigma)$ if and only
if $\Sigma$ is acyclic.
\end{theorem}

Combining Theorem~\ref{th:lower-bound-exact},
Theorem~\ref{th:acyclic-sufficient-lower},
and Corollary~\ref{cor:Grobner}, we obtain the following corollary.

\begin{corollary}
\label{cor:cluster-algebra-finite-generation}
Let $\AA = \AA(\Sigma)$ be the cluster algebra associated with
a totally mutable acyclic seed~$\Sigma$.
Then:
\begin{itemize}
\item $\AA$ is generated by $x_1, x'_1, \ldots, x_n, x'_n$.

\item The standard monomials in $x_1, x'_1, \ldots,
x_n, x'_n$ form a $\ZP$-basis of $\AA$.

\item The polynomials $x_j x'_j - P_j(\xx)$, for $j \in [1,n]$, generate the
ideal $\mathcal{I}$ of relations
among the generators $x_1,\dots,x_n,x_1',\dots,x_n'$.

\item These polynomials form a Gr\"obner basis
for~$\mathcal{I}$ with respect to any term order
in which $x_1',\dots,x_n'$ are much more expensive than $x_1,\dots,x_n\,$.
\end{itemize}
\end{corollary}

\begin{remark}
Corollary~\ref{cor:cluster-algebra-finite-generation} applies
to the cluster algebras of finite type classified in~\cite{fz-clust2}.
According to that classification, any such algebra $\AA$ has a seed
$\Sigma$ such that the corresponding graph $\Gamma(B)$ is an orientation of a
Dynkin diagram.
Since every Dynkin diagram is a tree, the seed $\Sigma$ is
acyclic, and Corollary~\ref{cor:cluster-algebra-finite-generation}
provides an explicit description of the cluster algebra~$\AA$.

If, in addition, the seed $\Sigma$ is coprime (which can always be achieved
by an appropriate choice of the coefficient tuple~$\pp$), then
$\AA = \overline\AA$ by Corollary~\ref{cor:finite-generation-sufficient}.
\end{remark}

\subsection{Open problems}

If a mutation class lacks an acyclic representative---and this does
happen a lot---then the structure of the corresponding cluster
algebra (or upper cluster algebra, or the lower/upper bounds)
is much less understood.
In this section, we discuss some of the natural questions for which we
only have partial answers.

\begin{problem}
Classify finitely generated (over $\ZP$) cluster algebras.

\end{problem}


This problem is trivial in rank $n\leq 2$ because all the seeds are acyclic.
In rank $n=3$, we obtain the following result.

 \begin{theorem}
\label{th:fin-gen=tot-cyclic}
A skew-symmetrizable cluster algebra of rank~$3$ is finitely generated
over~$\ZP$ if and only if it has an acyclic seed.
\end{theorem}

\begin{problem}
\label{problem:when-upper=cluster}
When is a cluster algebra equal to the corresponding upper cluster
algebra?
\end{problem}

The following statement shows that, already in rank~$3$,
the two algebras do not have to coincide.

\begin{proposition}
\label{pr:Conway}
Let $\Sigma$ be a seed with the (skew-symmetric) exchange matrix
\[
B=\left[\begin{array}{ccc}
0  &  2  & -2\\
-2 & 0 & 2\\
2  & -2  & 0 \\
\end{array}\right].
\]
Then $\AA(\Sigma)\neq\overline\AA(\Sigma)$.
\end{proposition}

\begin{problem}
\label{prob:finite-gen-upper}
When is an upper cluster algebra finitely generated?
\end{problem}

It follows from Corollaries~\ref{cor:finite-generation-sufficient}
and~\ref{cor:cluster-algebra-finite-generation}
that any upper cluster algebra that possesses a coprime and acyclic seed
is finitely generated.
However, the acyclicity condition is not necessary:
there exist finitely generated upper cluster algebras without an
acyclic seed
(see Remark~\ref{rem:finite-generation-double-cells}).
In fact, we do not know a single example of an upper cluster algebra
proved to be infinitely generated.

\begin{problem}
\label{prob:acyclic-seeds}
Is any acyclic seed totally mutable?
In other words, is any acyclic sign-skew-symmetric integer matrix
totally sign-skew-symmetric?
\end{problem}

A complete classification of totally sign-skew-symmetric integer
$3 \times 3$ matrices (thus, cluster algebras of rank~$3$)
obtained by the authors will appear in a separate publication.
In particular, this classification implies the affirmative answer to
Problem~\ref{prob:acyclic-seeds} for the cluster algebras of rank~$n \leq 3$.


\section{Cluster algebras associated with double Bruhat cells}
\label{sec:double-cells}

In this section, we relate the constructions of
Section~\ref{sec:main-results}
to the study of double Bruhat cells.
More concretely, we identify the coordinate rings of double Bruhat
cells with the upper
cluster algebras defined by appropriate combinatorial data.
The results in this section are furnished with complete proofs, which
draw heavily from the arguments and constructions in~\cite{fz-double},
\cite{ssvz}, and~\cite{z-imrn},
and rely on the results in Section~\ref{sec:main-results}.


\subsection{Double Bruhat cells}

To state the main result of Section~\ref{sec:double-cells}
(Theorem~\ref{double-cell=upper-bound-explicit} below),
we need to recall the relevant background
from~\cite{fz-double,z-imrn}.
Let $G$ be a simply connected, connected, semisimple complex
algebraic group of rank~$r$.
Let $B$ and $B_-$ be two opposite Borel subgroups in $G$, let $N$
and $N_-$ be their unipotent radicals, $H = B \cap B_-$ a maximal
torus, and $W = {\rm Norm}_G (H) / H$ the \emph{Weyl group}.
Let $\gg = {\rm Lie}(G)$ be the Lie algebra of $G$, and $\hh = {\rm Lie}(H)$
the Cartan subalgebra of $\gg$.
Let $\Pi = \{\alpha_1, \dots,
\alpha_r\} \subset \hh^*$ be the system of simple roots for which the
corresponding root subgroups are contained in~$N$.
Let $\alpha_i^\vee \in \hh$ denote the simple coroot corresponding to
a simple root $\alpha_i$, and let $A = (a_{ij})$
be the \emph{Cartan matrix} given by $a_{ij} = \alpha_j (\alpha_i^\vee)$.
For every $i \in [1,r]$,
let $\varphi_i: SL_2
\to G$ denote the corresponding 
$SL_2$-embedding that sends the upper- (resp., lower-) triangular matrices
in $SL_2$ to~$B$ (resp., to~$B_-$).

The Weyl group $W$ is canonically identified with the Coxeter
group generated by the involutions $s_1, \dots, s_r$
subject to the relations $(s_i s_j)^{d_{ij}}=e$ (here $e\in W$ is the
identity element) for all
$i\! \neq\!j$, where~$d_{ij} \!=\!2$ (resp., $3,4,6$) if $a_{ij}a_{ji} = 0$
(resp., $1,2,3$).
Under this identification, $s_i = \overline {s_i} H$, where
$$
\overline {s_i} = \varphi_i \mat{0}{-1}{1}{0} \in {\rm Norm}_G (H)\,.
$$

A word $\ii = (i_1, \ldots, i_\ell)$ in
the alphabet $[1,r]$ is a \emph{reduced word} for $w \in W$ if
$w = s_{i_1} \cdots s_{i_\ell}$, and $\ell$ is the smallest length of
such a factorization.
The length $\ell$ of any reduced word for $w$ is
called the \emph{length} of $w$ and denoted by $\ell (w)$.
We denote by $w_\circ$ the element of maximal length in~$W$.

The representatives $\overline {s_i} \in G$ satisfy the braid
relations
\[
\overline {s_i} \ \overline {s_j}\ \overline {s_i}
\cdots = \overline {s_j} \ \overline {s_i} \ \overline {s_j}
\cdots,
\]
with $d_{ij}$ factors on each side.
Thus, the
representative $\overline w \in {\rm Norm}_G (H)$ can be
unambiguously defined for any $w \in W$ by requiring that
$\overline {uv} = \overline {u} \cdot \overline {v}$ whenever
$\ell (uv) = \ell (u) + \ell (v)$.

In the special case of $G=SL_{r+1}(\CC)$ (type~$A_r$),
the Weyl group is naturally identified with the symmetric
group~$S_{r+1}$.
The permutation $w_\circ$ of maximal length (i.e., maximal number of
inversions) is given by $w_\circ(i)=r+2-i$.
The corresponding representative $\overline{w_\circ}$ is the
antidiagonal matrix
\[
\begin{bmatrix}
\cdots&\cdots&\cdots&\cdots\\
0 & 0 & 1&\cdots\\
0 & -1& 0&\cdots\\
1 & 0 & 0&\cdots\\
\end{bmatrix}\,.
\]

The group $G$ has two \emph{Bruhat decompositions}, with respect
to $B$ and $B_-\,$:
$$G = \bigcup_{u \in W} B u B = \bigcup_{v \in W} B_- v B_- \ . $$
The \emph{double Bruhat cells}~$G^{u,v}$ are
defined by $G^{u,v} = B u B \cap B_- v B_- \,$.
According to \cite[Theorem~1.1]{fz-double}, $G^{u,v}$ is
biregularly isomorphic
to a Zariski open subset of an affine space of dimension
$r+\ell(u)+\ell(v)$.

\begin{example}
\label{example:SL3-w0w0}
Our running example will be the open double Bruhat cell of
type~$A_2$.
More precisely, let $G=SL_3(\CC)$, let $B$ (resp.,~$B_-$) be the Borel
subgroup of upper- (resp., lower-) triangular matrices,
and let
$u=v=w_\circ=s_1s_2s_1=s_2s_1s_2$ be the order-reversing
permutation~$321$ (the element of maximal length in the symmetric
group~$W=S_3$).
Applying \cite[Proposition~4.1]{fz-double},
we see that the double Bruhat cell $G^{w_\circ,w_\circ}\subset SL_3(\CC)$
consists of all complex $3\times 3$ matrices $x=(x_{ij})$
of determinant~$1$ satisfying the conditions
\begin{equation}
\label{eq:sl3-w0w0-nonvanishing}
x_{13}\neq 0,
\quad
\Bigl|\!\Bigl|\begin{array}{cc}
x_{12} & x_{13} \\
x_{22} & x_{23}
\end{array}\Bigr|\!\Bigr|\neq 0,
\quad
x_{31}\neq 0,
\quad
\Bigl|\!\Bigl|\begin{array}{cc}
x_{21} & x_{22} \\
x_{31} & x_{32}
\end{array}\Bigr|\!\Bigr|\neq 0.
\end{equation}
\end{example}

\subsection{Combinatorial data}

The group $W \times W$ is a Coxeter group as well.
To avoid confusion, we will use the indices $-1, \dots, -r$ for the simple
reflections in the first copy of~$W$, and $1, \dots, r$ for the second copy.
Thus, a reduced word for a pair $(u,v) \in W \times
W$ is an arbitrary shuffle of a reduced word for $u$ written in
the alphabet $-[1,r]$ and a reduced word for $v$ written in the alphabet $[1,r]$.
For an index $i \in \pm[1,r]$, we denote by $\veps(i)$ the sign of~$i$.
We define the \emph{support} $\Supp(u,v) \subset [1,r]$
of a pair $(u,v)$ as the set of all indices $i$ such
that some (equivalently, any) reduced word for $(u,v)$ contains either
$i$ or~$-i$ (or both).

Let us fix a pair $(u,v) \in W \times W$ and a reduced word
$\ii = (i_1, \dots, i_{ \ell(u)+\ell(v)})$ for~$(u,v)$.
 From~$\ii$, we will now construct a rectangular matrix
$\tilde B = \tilde B(\ii)$ satisfying the conditions of
Corollary~\ref{cor:upper-bound-from-matrix}.
We will eventually identify the coordinate ring of $G^{u,v}$
with the corresponding upper cluster algebra $\overline\AA=\overline\AA(\tilde B)$
of geometric type.

We will construct the matrix $\tilde B(\ii)
$ in two steps.
First, in Definition~\ref{def:edges} we will build a directed graph
$\tilde \Gamma(\ii)$
that will determine the signs of the matrix entries of~$\tilde B(\ii)$.
We will then define the absolute values of these entries in
Definition~\ref{def:tildeB-double-cells}.
In the simply-laced cases (i.e., when every simple component of $G$ is
of type $A_r$, $D_r$, or~$E_r$),
all these absolute values will be equal to $0$ or~$1$,
so the matrix $\tilde B(\ii)$ will be completely determined by the
graph~$\tilde \Gamma(\ii)$.

Let us add $r$ additional entries $i_{-r}, \ldots, i_{-1}$ at the
beginning of $\ii$ by setting $i_{-j} = -j$ for $j \in [1,r]$
(note that the index $0$ is missing).
The rows of the matrix $\tilde B(\ii)$ will be labeled by the
elements of the set
$-[1,r] \cup [1,\ell(u)+\ell(v)]$,
which has cardinality
\[
m = r+\ell(u)+\ell(v).
\]
For $k \in -[1,r] \cup [1,\ell(u)+\ell(v)]$, we denote by $k^+ =
k^+_\ii$ the smallest index $\ell$ such that $k < \ell$ and $|i_\ell|
= |i_k|$;
if $|i_k| \neq |i_\ell|$ for $k < \ell$, then we set $k^+ =
\ell(u)+\ell(v)+1$.
We say that an index $k$ is \emph{$\ii$-exchangeable} if both $k$ and $k^+$ belong to
$[1,\ell(u)+\ell(v)]$.
The columns of $\tilde B(\ii)$ will be labeled by the
elements of the set
$\ee (\ii)$ of all $\ii$-exchangeable indices; its cardinality is
\[
n =\ell(u)+\ell(v) - |\Supp(u,v)|.
\]
Thus, the matrix $\tilde B = \tilde B(\ii)$ will be of size $m\times
n$.
(By convention, if $n=0$, that is,
$|i_k| \neq |i_\ell|$ for $k \neq \ell$, then $\overline\AA$
is simply the ring of Laurent polynomials in $m$ variables.)

The following definition is a slight modification of the constructions
in \cite{ssvz,z-imrn}.

\begin{definition}
\label{def:edges}
We associate with $\ii$ a directed graph $\tilde \Gamma(\ii)$.
The set of vertices of $\tilde \Gamma(\ii)$ is $ -[1,r] \cup [1,\ell(u)+\ell(v)]$.
Vertices $k$ and~$\ell$, with $k<\ell$, are connected by an edge in
$\tilde \Gamma(\ii)$ if and only if either $k$ or~$\ell$ (or both) are
$\ii$-exchangeable, and one of the following three
conditions is satisfied:

\begin{enumerate}

\item
$\ell=k^+$;

\item
$\ell < k^+ < \ell^+$, $a_{|i_k|, |i_\ell|} < 0$, and
$\veps(i_{\ell})=\veps(i_{k^+})$;

\item
$\ell < \ell^+ < k^+$, $a_{|i_k|, |i_\ell|} < 0$, and
$\veps(i_{\ell})=-\veps(i_{\ell^+})$.

\end{enumerate}
The edges of type (1) are called \emph{horizontal}, and those of
types (2) and (3) are called \emph{inclined}.
A horizontal (resp., inclined)
edge is directed from $k$ to $\ell$ (as above, $k<\ell$)
if and only if $\veps (i_{\ell}) = +1$ (resp., $\veps (i_{\ell}) = -1$).
\end{definition}

\begin{definition}
\label{def:tildeB-double-cells}
We define an integer $m\times n$ matrix $\tilde B = \tilde B(\ii)$ with
rows labeled by all the indices in $-[1,r] \cup [1,\ell(u)+\ell(v)]$ and columns
labeled by the $\ii$-exchangeable indices in~$\ee(\ii)$.
A matrix entry $b_{k \ell}$ of $\tilde B$ is determined by the following rules:

\begin{enumerate}
\item
$b_{k\ell} \neq 0$ if and only if there is an edge
of $\tilde \Gamma(\ii)$ connecting $k$ and~$\ell$; \\
$b_{k\ell} > 0$ (resp., $b_{k\ell}<0$) if this edge
is directed from $k$ to~$\ell$ (resp., from $\ell$ to~$k$);

\item If $k$ and $\ell$ are connected by an edge of $\tilde \Gamma(\ii)$, then
\begin{equation}
\label{eq:absvalue-bkl}
|b_{k \ell}| =
\begin{cases}
1 & \text{if $|i_k|=|i_\ell|$\quad (a horizontal edge);} \\[.05in]
-a_{|i_k|,|i_\ell|} & \text{if $|i_k|\neq|i_\ell|$\quad (an inclined edge).}
\end{cases}
\end{equation}
\end{enumerate}
\end{definition}

\begin{remark}
\label{rem:Btilde-direct}
The matrix $\tilde B(\ii)$ can be computed directly from $\ii$ as
follows.
Let $k\in -[1,r] \cup [1,\ell(u)+\ell(v)]$ and $\ell\in\ee(\ii)$.
Set $p=\max(k,l)$ and $q=\min(k^+,\ell^+)$.
Then
\[
b_{k\ell} =
\begin{cases}
-\sgn(k-\ell)\cdot\varepsilon(i_p) & \text{if $p=q$;}\\
-\sgn(k-\ell)\cdot\varepsilon(i_p)\cdot a_{|i_k|,|i_\ell|}
   & \text{if $p<q$ and $\varepsilon(i_p) \varepsilon(i_q) (k-\ell)(k^+-\ell^+)>0$;}\\
0 & \text{otherwise.}
\end{cases}
\]
\end{remark}

\begin{example}
\label{example:SL3-w0w0-continued}
We continue with Example~\ref{example:SL3-w0w0}.
Thus, $r=2$, $u=v=w_\circ\in S_3$, $\ell(u)=\ell(v)=3$, $m=8$, and $n=4$.
Let us take $\ii=(1,2,1,-1,-2,-1)$.
Then $\ee(\ii)=[1,4]$.
The resulting graph $\tilde \Gamma(\ii)$ is shown in
Figure~\ref{fig:gamma-sl3}, where the vertices
$-2,-1,1,2,3,4,5,6$ can be traced from left to right (in this order),
each vertex $k$ is placed at height~$|i_k|$,
and the sign of~$i_k$ is indicated at each vertex.

The $8\times 4$ matrix $\tilde B(\ii)$, with its rows and columns
labeled by the sets $[-2,-1]\cup [1,6]$ and $[1,4]$, is shown in
Figure~\ref{fig:tildeB-sl3} on the left.
\end{example}

\begin{figure}[ht]
\setlength{\unitlength}{3pt}
\begin{picture}(70,15)(0,-3)
\put(2,10){\line(1,0){26}}
\put(32,10){\line(1,0){26}}
\put(12,0){\line(1,0){6}}
\put(22,0){\line(1,0){16}}
\put(42,0){\line(1,0){6}}
\put(52,0){\line(1,0){16}}
\put(2,9){\line(2,-1){16}}
\put(32,9){\line(2,-1){16}}
\put(21.5,1.5){\line(1,1){7}}
\put(51.5,1.5){\line(1,1){7}}
\put(38.5,1.5){\line(-1,1){7}}

\put(15,10){\vector(1,0){1}}
\put(45,10){\vector(-1,0){1}}
\put(15,0){\vector(1,0){1}}
\put(30,0){\vector(1,0){1}}
\put(45,0){\vector(-1,0){1}}
\put(60,0){\vector(-1,0){1}}
\put(10,5){\vector(-2,1){1}}
\put(40,5){\vector(2,-1){1}}
\put(25,5){\vector(-1,-1){1}}
\put(35,5){\vector(-1,1){1}}
\put(55,5){\vector(1,1){1}}

\put(0,10){\circle{4}}
\put(0,10){\makebox(0,0){$\scriptstyle -$}}
\put(30,10){\circle{4}}
\put(30,10){\makebox(0,0){$\scriptstyle +$}}
\put(60,10){\circle{4}}
\put(60,10){\makebox(0,0){$\scriptstyle -$}}
\put(10,0){\circle{4}}
\put(10,0){\makebox(0,0){$\scriptstyle -$}}
\put(20,0){\circle{4}}
\put(20,0){\makebox(0,0){$\scriptstyle +$}}
\put(40,0){\circle{4}}
\put(40,0){\makebox(0,0){$\scriptstyle +$}}
\put(50,0){\circle{4}}
\put(50,0){\makebox(0,0){$\scriptstyle -$}}
\put(70,0){\circle{4}}
\put(70,0){\makebox(0,0){$\scriptstyle -$}}

\end{picture}
\caption{Graph $\tilde \Gamma(\ii)$ for $\ii=(1,2,1,-1,-2,-1)$}
\label{fig:gamma-sl3}
\end{figure}

\begin{figure}[ht]
\[
\begin{array}{c|rrrr}
   &  1 &  2 &  3 &  4 \\
\hline
-2 & -1 &  1 &  0 &  0\\
-1 &  1 &  0 &  0 &  0\\
 1 &  0 & -1 &  1 &  0\\
 2 &  1 &  0 & -1 &  1\\
 3 & -1 &  1 &  0 & -1\\
 4 &  0 & -1 &  1 &  0\\
 5 &  0 &  1 &  0 & -1\\
 6 &  0 &  0 &  0 &  1
\end{array}
\qquad\qquad
\begin{array}{c|rrrr}
   &  1 &  2 &  3 &  4 \\
\hline
-1 &  1 &  0 &  0 &  0 \\
-2 & -1 &  1 &  0 &  0 \\
 1 &  0 & -1 &  1 &  0 \\
 3 & -1 &  1 &  0 & -1
\end{array}
\]
\caption{Matrix $\tilde B(\ii)$
and its $\ee(\ii)^-\!\times\ee(\ii)$ submatrix
for $\ii=(1,2,1,-1,-2,-1)$}
\label{fig:tildeB-sl3}
\end{figure}

The \emph{principal part} of $\tilde B$ is the $n \times n$
submatrix $B = B(\ii)$ of $\tilde B$ formed by the rows and
columns labeled by the $\ii$-exchangeable indices in~$\ee (\ii)$.
By the rule~(1) of Definition~\ref{def:tildeB-double-cells},
the matrix $B$ is sign-skew-symmetric.
The corresponding graph $\Gamma (B)$ (see Section~\ref{sec:main-results})
is the induced subgraph of $\tilde \Gamma(\ii)$ on the set of vertices~$\ee(\ii)$.

The following result can be derived from \cite[Proposition~3.2]{ssvz}
and its proof given there. We present a self-contained proof below.

\begin{proposition}
\label{pr:double-cell-skewsymm-coprime}
The matrix $\tilde B = \tilde B(\ii)$ has full rank~$n$.
Its principal part $B = B(\ii)$ is skew-symmetrizable.
\end{proposition}

\begin{proof}
To prove the first statement, it suffices to show that the $n \times n$ minor of
$\tilde B$ with the row set $$\ee (\ii)^- = \{k \in -[1,r] \cup
[1,\ell(u)+\ell(v)] \, : \, k^+ \in \ee(\ii)\}$$ is nonzero.
In view of Definitions~\ref{def:edges}
and~\ref{def:tildeB-double-cells}, if $k \in \ee (\ii)^-$ and
$\ell \in \ee (\ii)$, then $|b_{k \ell}| =1$ if $k^+ = \ell$,
and $|b_{k \ell}| =0$ if $k^+ < \ell$,
implying the claim.
(Cf.\ Figure~\ref{fig:tildeB-sl3}.
Note that we swapped the rows labeled by $-1$ and~$-2$,
to make the ordering of rows and columns compatible with the
bijection $k\mapsto k^+$.)

To prove the second statement, recall that the Cartan matrix $A$ is
symmetrizable, i.e., $d_i a_{ij} = d_j a_{ji}$ for some positive
numbers $d_1, \dots, d_r$.
In view of (\ref{eq:absvalue-bkl}), we
have $d_{|i_k|} |b_{k \ell}| = d_{|i_\ell|} |b_{\ell k}|$ for all
$k$ and $\ell$, which implies that $B$ is skew-symmetrizable.
\end{proof}

Proposition~\ref{pr:double-cell-skewsymm-coprime} shows that
Corollary~\ref{cor:upper-bound-from-matrix} applies to the
case under consideration.
Thus, the matrix $\tilde B(\ii)$
gives rise to a well-defined upper cluster
algebra $\overline\AA(\ii)$ of geometric type that coincides with the
upper bound~$\Upper(\Sigma)$ associated with a seed
$\Sigma(\ii)=(\xx,\tilde B(\ii))$ of geometric type.
Here we need to adjust the setup described in Section~\ref{sec:main-results}
a little bit, so that the cluster variables in $\xx$
are labeled by the index set~$\ee(\ii)$, and the
generators of the coefficient group~$\PP$ by the remaining indices
in
$-[1,r]\cup [1,\ell(u)+\ell(v)] - \ee (\ii)$.
Accordingly, the ambient field
$\FFcal$ of $\overline\AA(\ii)$ is the field of rational functions
(over~$\QQ$) in the $m$ independent variables
\begin{equation}
\label{eq:tilde-x}
\tilde \xx = \{x_k \, : \, k \in
-[1,r] \cup [1,\ell(u)+\ell(v)]\}\,.
\end{equation}

\begin{example}
In Example~\ref{example:SL3-w0w0-continued},
$\PP$ is the multiplicative group generated by the variables
$x_{-2},x_{-1},x_5,x_6$, whereas the cluster~$\xx$ associated
with~$\ii$ consists of $x_1,x_2,x_3,x_4$.
The exchange relations~\eqref{eq:exchange} are produced by looking at
the columns of the matrix~$\tilde B(\ii)$ in
Figure~\ref{fig:tildeB-sl3}.
They are:
\begin{equation}
\label{eq:exchange-sl3w0w0}
\begin{array}{rcl}
x_1 x_1'&=&x_{-1}x_2+x_{-2}x_3\,,\\
x_2 x_2'&=&x_{-2}x_3x_5+x_1x_4\,,\\
x_3 x_3'&=&x_1x_4+x_2\,,\\
x_4 x_4'&=&x_2x_6+x_3x_5\,.
\end{array}
\end{equation}
The 
algebra $\overline\AA(\ii)$ consists of all
rational functions in $\FFcal=\QQ(x_{-2},x_{-1},x_1,\dots,x_6)$
that can be written as Laurent polynomials in each of the five
clusters
\[
\xx=(x_1,x_2,x_3,x_4), \ \xx_1=(x_1',x_2,x_3,x_4), \dots,
\ \xx_4=(x_1,x_2,x_3,x_4').
\]
\end{example}

\subsection{Generalized minors}

Let $\overline\AA(\ii)_\CC =
\overline\AA(\ii) \otimes \CC$ be the $\CC$-algebra obtained
from $\overline\AA(\ii)$ by extension of scalars.
Its ambient field is $\FFcal_\CC=\FFcal\otimes\CC$.
Our goal is to show that
the coordinate ring $\CC[G^{u,v}]$ is naturally isomorphic
to~$\overline\AA(\ii)_\CC$.
To describe this isomorphism, we will need to recall from
\cite{fz-double} some terminology and background concerning
\emph{generalized minors}.
First of all, the \emph{weight lattice}
$P$ of $G$ can be thought of as the group of rational
multiplicative characters of~$H$ written in the exponential
notation: a weight $\gamma\in P$ acts by $a \mapsto a^\gamma$.
The lattice $P$ is also identified with the additive group of all
$\gamma \in \hh^*$ such that $\gamma (\alpha_i^\vee) \in \ZZ$ for
all $i \in [1,r]$.
Thus, $P$ has a $\ZZ$-basis $\{\omega_1, \dots,
\omega_r\}$ of \emph{fundamental weights} given by $\omega_j
(\alpha_i^\vee) = \delta_{i,j}$.

Denote by $G_0=N_- H N$ the open subset of elements $x\in G$ that
have Gaussian decomposition; this (unique) decomposition will be
written as $x = [x]_- [x]_0 [x]_+ \,$, where $[x]_- \in N_-$,
$[x]_0 \in H$, and $[x]_+ \in N$.
For $u,v \in W$ and $i \in
[1,r]$, the generalized minor $\Delta_{u \omega_i, v \omega_i}$
is the regular function on $G$ whose restriction to the open set
${\overline {u}} G_0 {\overline {v}}^{-1}$ is given by
\begin{equation}
\label{eq:Delta-general} \Delta_{u \omega_i, v \omega_i} (x) =
(\left[{\overline {u}}^{\ -1}
  x \overline v\right]_0)^{\omega_i} \ .
\end{equation}
As shown in \cite{fz-double}, $\Delta_{u \omega_i, v \omega_i}$
depends on the weights $u \omega_i$ and $v \omega_i$ alone, not on
the particular choice of $u$ and~$v$.
It is easy to see that the
generalized minors are distinct irreducible elements of the
coordinate ring $\CC[G]$.

In the type~$A_r$ special case,
where $G=SL_{r+1}(\CC)$ is the group of $(r+1)\times(r+1)$
complex matrices of determinant~$1$ (with the standard choices of
opposite Borel subgroups),
the generalized minors are nothing but the ordinary minors of a matrix.
More precisely, for permutations $u,v\in S_{r+1}$, the evaluation of
$\Delta_{u \omega_i, v \omega_i}$ at a matrix $x\in SL_{r+1}$
is equal to the determinant of the submatrix of $x$ whose rows
(resp., columns) are labeled by the elements of the set
$u([1,i])$ (resp.,~$v([1,i])$).

Each Bruhat cell, and therefore each double Bruhat cell,
can be defined inside~$G$ by a collection of vanishing/non-vanishing conditions
of the form $\Delta(x)=0$ and \hbox{$\Delta(x)\neq0$,}
where $\Delta$ is a generalized minor.
A precise statement appears in Proposition~\ref{pr:double-cell-by-eqs}
below. We follow~\cite[Section~3]{fz-recogn},
although the result itself goes back to~\cite{GP1,GP5}.
We refer the reader to \cite{fz-recogn} and
\cite[Section~2]{fz-double} for further details.

For a fixed $i\in[1,r]$, the weights of the form $w\omega_i$
are in a canonical bijection with the cosets in $W$
modulo the stabilizer~$W_{\hat i}$ of~$\omega_i$.
Explicitly, $W_{\hat i}$ is the maximal parabolic subgroup
generated by the simple reflections~$s_j$ with~$j\neq i$.
The restriction of the Bruhat order on $W$ to the set of minimal coset
representatives (modulo~$W_{\hat i}$) induces a partial order on the
weights~$w\omega_i$, which is also called the Bruhat order. \linebreak[2]
In the special case $G=SL_{r+1}(\CC)$, this Bruhat order translates into
the partial order on the $i$-element subsets of~$[1,r+1]$
defined by setting
\[
\{j_1<\cdots<j_i\}\leq\{k_1<\cdots<k_i\}
\stackrel{\rm def}{\Longleftrightarrow}
(j_1\leq k_1,\dots,j_i\leq k_i)\,.
\]

\begin{proposition}
\label{pr:double-cell-by-eqs}
A double Bruhat cell $G^{u,v}$ is given inside~$G$ by the following
conditions (here $i$ runs over the set~$[1,r]$):
\begin{align}
& \Delta_{u'\omega_i,\omega_i}=0 \quad \text{whenever $u'\omega_i\not\leq
  u\omega_i$ in the Bruhat order,}\\
& \Delta_{\omega_i,v'\omega_i}=0 \quad \text{whenever $v'\omega_i\not\leq
  v^{-1}\omega_i$ in the Bruhat order,}\\
\label{eq:double-cell-ineq}
&\Delta_{u\omega_i,\omega_i}\neq 0, \quad
  \Delta_{\omega_i,v^{-1}\omega_i}\neq 0\,.
\end{align}
\end{proposition}

To illustrate, the non-vanishing conditions \eqref{eq:sl3-w0w0-nonvanishing}
are a special case of~\eqref{eq:double-cell-ineq}.

More economical descriptions of the double Bruhat cells can be
obtained from \cite[Proposition~4.1]{fz-double}
and~\cite[Section~4]{fz-recogn}.

Slightly modifying the construction in
\cite[Section~1.6]{fz-double}, we introduce the following notation.
For $k \in [1,\ell(u)+\ell(v)]$, let us denote
\begin{align}
\label{eq:u-partial}
& u_{\leq k} = u_{\leq k}(\ii) =
\doublesubscript{\prod}{\ell = 1, \ldots, k} {\veps(i_\ell) = -1}
s_{|i_\ell|} \, , \\[.1in]
\label{eq:v-partial}
& v_{>k} = v_{>k}(\ii) =
\doublesubscript{\prod}{\ell = \ell(u)+\ell(v), \ldots, k+1}
{\veps(i_\ell) = +1} s_{|i_\ell|} \,,
\end{align}
where in the noncommutative product in (\ref{eq:u-partial}) (resp.,
(\ref{eq:v-partial})),
the index $\ell$ is increasing (resp., decreasing).
For $k \in -[1,r]$, we adopt the convention that $u_{\leq k}=e$ and
$v_{>k}=v^{-1}$.
Now for $k \in -[1,r] \cup [1,\ell(u)+\ell(v)]$,
we set
\begin{equation}
\label{eq:Delta-factors}
\Delta(k;\ii) = \Delta_{u_{\leq k}
\omega_{|i_k|}, v_{>k} \omega_{|i_k|}} \, .
\end{equation}
Finally, we denote by
\begin{equation}
\label{eq:F(i)}
F(\ii) = \{ \Delta (k;\ii) : k \in -[1,r] \cup [1,\ell(u)+\ell(v)]\}
\end{equation}
the entire set of
$r+\ell(u)+\ell(v)$ minors associated with a fixed reduced
word~$\ii$.

A word of warning: the notation \eqref{eq:Delta-factors}
is different from the one used in
\cite{fz-double}, where the minors $\Delta (k;\ii)$ (in our
present notation) were associated with the reduced word for
$(u^{-1}, v^{-1})$ obtained by reading $\ii$ backwards.
On the other hand, \eqref{eq:F(i)}
is consistent with~\cite[(1.22)]{fz-double}.

\begin{example}
\label{example:minors-sl3-w0w0}
Figure~\ref{fig:minors-sl3} illustrates the definitions
\eqref{eq:u-partial}--\eqref{eq:Delta-factors}
for our running example.
We represent each weight $u_{\leq k}\omega_{|i_k|}$
(resp., $v_{>k} \omega_{|i_k|}$) by the corresponding subset of rows
(resp., columns) $u_{\leq k}([1,k])$ (resp.,~$v_{>k}([1,k])$).
Note that, in this example, the collection of minors $F(\ii)$
consists, in the terminology of~\cite{fz-intel},
of all ``initial'' minors of a $3\times 3$ matrix
$x=(x_{ij})$, except for $\det(x)=1$.
These are the minors that occupy consecutive rows and columns,
including at least one entry in the first row or the first column.
\end{example}

\begin{figure}[ht]
\[
\begin{array}{c|cccccccc}
k&-2&-1&1&2&3&4&5&6\\
\hline
i_k&-2&-1&1&2&1&-1&-2&-1\\
\hline
u_{\leq k} & e & e & e & e & e & s_1 & s_1s_2 & w_\circ\\
v_{>k} &w_\circ & w_\circ & s_1s_2 & s_1 & e & e & e & e\\
\hline
u_{\leq k}\omega_{|i_k|} & 12 & 1&1&12&1&2&23&3\\
v_{>k} \omega_{|i_k|} & 23&3&2&12&1&1&12&1\\
\hline
\\[-.15in]
\Delta(k;\ii) & \biggl|\!\biggl|\!\begin{array}{cc}
x_{12} & x_{13} \\
x_{22} & x_{23}
\end{array}\!\biggr|\!\biggr|
& x_{13}
& x_{12}
&\biggl|\!\biggl|\!\begin{array}{cc}
x_{11} & x_{12} \\
x_{21} & x_{22}
\end{array}\!\biggr|\!\biggr|
& x_{11}
& x_{21}
& \biggl|\!\biggl|\!\begin{array}{cc}
x_{21} & x_{22} \\
x_{31} & x_{32}
\end{array}\!\biggr|\!\biggr|
& x_{31}
\end{array}
\]
\caption{Minors $\Delta(k;\ii)$
for $\ii=(1,2,1,-1,-2,-1)$}
\label{fig:minors-sl3}
\end{figure}

By~\cite[Theorem~1.12]{fz-double}, the
set $F(\ii)$ is a transcendence basis for the field of
rational functions~$\CC(G^{u,v})$.


\subsection{Cluster algebra structures in $\CC[G^{u,v}]$}
We can now state the main result of Section~\ref{sec:double-cells}.

\begin{theorem}
\label{double-cell=upper-bound-explicit}
For every reduced word
$\ii$ of a pair $(u,v)$ of Weyl group elements, the isomorphism of
fields $\varphi: \FFcal_\CC \to\CC(G^{u,v})$ defined by
\[
\text{$\varphi(x_k) = \Delta(k,\ii)$\ \ for all\ \ $k \in -[1,r] \cup
[1,\ell(u)+\ell(v)]$}
\]
restricts to an isomorphism of algebras
$\overline\AA(\ii)_\CC \to \CC[G^{u,v}]$.
\end{theorem}

\begin{example}
We continue with Example~\ref{example:minors-sl3-w0w0}.
In this example,
the isomorphism~$\varphi$ of Theorem~\ref{double-cell=upper-bound-explicit}
identifies the coefficient group~$\PP$ with the group generated by the
minors $\Delta(k,\ii)$ with $k\notin\ee(\ii)$, i.e., $k\in\{-2,-1,5,6\}$.
Thus, $\varphi(\PP)$ consists of all Laurent monomials in the minors
\begin{equation}
\label{eq:4-frozen-minors}
\biggl|\!\biggl|\!\begin{array}{cc}
x_{12} & x_{13} \\
x_{22} & x_{23}
\end{array}\!\biggr|\!\biggr|,\ \
x_{13},\ \
\biggl|\!\biggl|\!\begin{array}{cc}
x_{21} & x_{22} \\
x_{31} & x_{32}
\end{array}\!\biggr|\!\biggr|,\ \
x_{31}.
\end{equation}
Note that by~\eqref{eq:sl3-w0w0-nonvanishing}, these are precisely the
minors that vanish nowhere on the double Bruhat cell under
consideration.

The cluster variables $x_1,x_2,x_3,x_4$
in the seed of geometric type $\Sigma(\ii)$
are identified by the isomorphism~$\varphi$ with the minors
\begin{equation}
\label{eq:4-minors}
x_{12},\ \
\biggl|\!\biggl|\!\begin{array}{cc}
x_{11} & x_{12} \\
x_{21} & x_{22}
\end{array}\!\biggr|\!\biggr|,\ \
x_{11}\,,\ \
x_{21}\,,
\end{equation}
respectively.
We next compute the images under~$\varphi$ of
the cluster variables $x_1',\dots,x_4'$,
using the exchange relations~\eqref{eq:exchange-sl3w0w0}.
For instance,
\begin{align}
\label{eq:varphi(x1')}
\varphi(x_1')&=
  \varphi\Bigl(\displaystyle\frac{x_{-1}x_2+x_{-2}x_3}{x_1}\Bigr)
\\
\nonumber
  &=\Bigl(x_{13}\biggl|\!\biggl|\!\begin{array}{cc}
x_{11} & x_{12} \\
x_{21} & x_{22}
\end{array}\!\biggr|\!\biggr|+\biggl|\!\biggl|\!\begin{array}{cc}
x_{12} & x_{13} \\
x_{22} & x_{23}
\end{array}\!\biggr|\!\biggr|x_{11}\Bigr)x_{12}^{-1} \\
\nonumber
&=
\biggl|\!\biggl|\!\begin{array}{cc}
x_{11} & x_{13} \\
x_{21} & x_{23}
\end{array}\!\biggr|\!\biggr|\,.
\end{align}
Similar calculations yield
\begin{align}
\label{eq:varphi(x2')}
\varphi(x_2') &=
x_{12}x_{21}x_{33}-x_{12}x_{23}x_{31}-x_{13}x_{21}x_{32}+x_{13}x_{22}x_{31}\,,\\
\label{eq:varphi(x3')}
\varphi(x_3') &= x_{22}\,,\\
\label{eq:varphi(x4')}
\varphi(x_4') &= \biggl|\!\biggl|\!\begin{array}{cc}
x_{11} & x_{12} \\
x_{31} & x_{32}
\end{array}\!\biggr|\!\biggr|\,.
\end{align}
(In order to obtain~\eqref{eq:varphi(x2')},
we need to use that $\det(x)=1$.)
Note that each $\varphi(x_i')$ is a regular function on
$G=SL_3(\CC)$.

In this particular example,
Theorem~\ref{double-cell=upper-bound-explicit} provides the following
description of the coordinate ring $\CC[G^{w_\circ,w_\circ}]$
of the open double Bruhat cell
$G^{w_\circ,w_\circ}\subset SL_3(\CC)$.
This coordinate ring consists precisely of those rational functions in the matrix
elements~$x_{ij}$ that can be written as Laurent polynomials,
with complex coefficients, in each of the following five transcendence
bases of the field~$\CC(x_{ij})$. The first basis is
\[
F(\ii)=\Biggl\{
\biggl|\!\biggl|\!\begin{array}{cc}
x_{12} & x_{13} \\
x_{22} & x_{23}
\end{array}\!\biggr|\!\biggr|,{\ }
x_{13}\,,{\ }
x_{12}\,,{\ }
\biggl|\!\biggl|\!\begin{array}{cc}
x_{11} & x_{12} \\
x_{21} & x_{22}
\end{array}\!\biggr|\!\biggr|,{\ }
x_{11}\,,{\ }
x_{21}\,,{\ }
\biggl|\!\biggl|\!\begin{array}{cc}
x_{21} & x_{22} \\
x_{31} & x_{32}
\end{array}\!\biggr|\!\biggr|,{\ }
x_{31}
\Biggr\}.
\]
The other four bases are obtained by replacing each of the minors
$\varphi(x_i)$ in~\eqref{eq:4-minors} by its counterpart $\varphi(x_i')$
as calculated in \eqref{eq:varphi(x1')}--\eqref{eq:varphi(x4')}.
\end{example}

\noindent
\textbf{Proof of Theorem~\ref{double-cell=upper-bound-explicit}.}
As in Section~\ref{sec:main-results}, we denote by
$x'_\ell$ the cluster variable obtained by exchanging $x_\ell$
from the initial seed in accordance with (\ref{eq:exchange});
the only difference is that now the index $\ell$ runs over the set $\ee(\ii)$.
We will need the following lemma.

\begin{lemma}
\label{lem:Delta-cluster-properties}
{\ }
\begin{enumerate}
\item
The minors $\Delta(k;\ii)$, for
$k\notin\ee (\ii)$,
vanish nowhere on~$G^{u,v}$.

\item
The map $G^{u,v} \to \CC^{r+\ell(u)+\ell(v)}$ defined by
$g\mapsto(\Delta(g))_{\Delta \in F(\ii)}$ restricts to a
biregular isomorphism $U(\ii) \to \CC_{\neq
0}^{r+\ell(u)+\ell(v)}$, where
\[
U(\ii)\stackrel{\rm def}{=}\{g \in G^{u,v} : \text{$\Delta(g) \neq 0$ for all $\Delta \in
 F(\ii)$}\}.
\]

\item
For each $\ell \in \ee(\ii)$, the rational function
\[
\Delta'(\ell;\ii) \stackrel{\rm def}{=} \varphi (x'_\ell) \in \CC(G^{u,v})
\]
is regular on $G^{u,v}$, i.e., it belongs to $\CC[G^{u,v}]$.

\item
For each $\ell \in \ee(\ii)$, the map $G^{u,v} \to
\CC^{r+\ell(u)+\ell(v)}$ defined by
$g\mapsto (\Delta(g))_{\Delta \in F_\ell (\ii)}$, where
\[
F_\ell(\ii) \stackrel{\rm def}{=} F(\ii) - \{\Delta(\ell;\ii)\} \cup
\{\Delta'(\ell;\ii)\},
\]
restricts to a biregular
isomorphism $U_\ell(\ii) \to \CC_{\neq 0}^{r+\ell(u)+\ell(v)}$,
where
\[
U_\ell(\ii)\stackrel{\rm def}{=}\{g \in G^{u,v}:
\text{$\Delta(g) \neq 0$ for all $\Delta \in F_\ell(\ii)$}\}.
\]
\end{enumerate}
\end{lemma}

\begin{proof}
This lemma is a direct consequence of \cite[Lemma~3.1]{z-imrn}
and its proof given there.
(In fact, statements (1) and (2) follow
directly from the results in~\cite{fz-double}.)
One only needs to adjust the setup in \cite{z-imrn} as follows.
First, in \cite{z-imrn}, the double cell $G^{u,v}$ was replaced by the
\emph{reduced} double cell (see~\cite[Proposition~4.3]{bz-invent})
$$L^{u,v} = \{g \in G^{u,v} \, : \, \Delta_{u\omega_i,
\omega_i}(g) = 1 \quad (i \in [1,r])\}.$$
Second, the results in
\cite{z-imrn} were stated in terms of ``twisted" generalized
minors obtained from ordinary minors by composing them with a certain
biregular isomorphism (the ``twist") between two (reduced) double
cells. We do not need this twist now because all elements of our
initial cluster are in fact actual (not twisted) generalized minors.
Translating the results of \cite{z-imrn} into our present
setup is straightforward, with most of the effort going into
proving part~(3).
This statement is a consequence of 
\cite[Theorem~4.3]{z-imrn}.
Recall that the latter theorem asserts that certain exchange-like
relations involving generalized minors produce regular functions on the entire group~$G$.
Now a direct check shows that the exchange relations that produce
the functions
$\Delta'(\ell;\ii)$ are obtained from the relations in \cite[Theorem~4.3]{z-imrn}
by applying them to an appropriate double reduced word (essentially,
the word~$\ii$ read backwards) and using the transformations of
minors in~\cite[Proposition~2.7]{fz-double}.
\end{proof}

The same argument as in \cite[Lemma~3.6]{z-imrn} shows that
Lemma~\ref{lem:Delta-cluster-properties} has the following
important corollary.
Let
\[
U = U(\ii) \ \cup \ \bigcup_{\ell \in \ee(\ii)} U_\ell(\ii).
\]
By definition, $U$ is a (Zariski) open subset in~$G^{u,v}$.

\begin{lemma}
\label{lem:complement}
The complement $G^{u,v} - U$ has (complex) codimension $\geq 2$ in $G^{u,v}$.
\end{lemma}

Now everything is ready for the proof of
Theorem~\ref{double-cell=upper-bound-explicit}.
Let us set $\tilde \xx_\ell = \tilde \xx -
\{x_\ell\} \cup \{x'_\ell\}$ (see \eqref{eq:tilde-x}).
We already observed that Corollary~\ref{cor:upper-bound-from-matrix}
applies to the upper bound $\overline\AA(\ii)$,
implying that
\begin{equation*}
\overline\AA(\ii)_\CC
=\CC[\tilde \xx^{\pm 1}] \cap \bigcap_{\ell \in \ee(\ii)} \CC
[\tilde \xx_\ell^{\pm 1}].
\end{equation*}
Hence it is enough to show that
\begin{equation}
\label{eq:upper-bound-concrete-2} \CC[G^{u,v}] =
\varphi(\CC[\tilde \xx^{\pm 1}]) \cap \bigcap_{\ell \in
\ee(\ii)} \varphi (\CC [\tilde \xx_\ell^{\pm 1}]).
\end{equation}
By Lemma~\ref{lem:Delta-cluster-properties}, Parts~2 and~4, we
have $$\varphi(\CC[\tilde \xx^{\pm 1}]) = \CC[U(\ii)], \quad
\varphi(\CC[\tilde \xx_\ell^{\pm 1}]) = \CC[U_\ell(\ii)],$$
so the right-hand side of (\ref{eq:upper-bound-concrete-2}) consists
of rational functions on $G^{u,v}$ that are regular on~$U$.
By Lemma~\ref{lem:complement}, any function regular on $U$ is regular
on the whole double cell~$G^{u,v}$, and we are done.
\endproof

\begin{remark}
\label{rem:ii-dependence}
Theorem~\ref{double-cell=upper-bound-explicit}
identifies the coordinate ring $\CC[G^{u,v}]$ with the upper bound
$\overline\AA(\ii)_\CC$ defined in terms of a reduced word $\ii$ for $(u,v)$.
More precisely, by Theorem~\ref{double-cell=upper-bound-explicit},
the initial seed of the cluster algebra $\AA(\ii)$ is determined
by the family of $m$ minors $F(\ii)$ and the $m \times n$ matrix
$\tilde B(\ii)$.
This raises the following natural question: if $\ii'$ is another
reduced word for $(u,v)$, is it true that the pair $(F(\ii'), \tilde B(\ii'))$
is 
a seed of $\AA(\ii)$, i.e., can be obtained
from $(F(\ii), \tilde B(\ii))$ by a sequence of mutations?
If $G$ is simply-laced, the positive answer to this
question follows from the results in~\cite{ssvz}.
To be more precise, in the simply-laced case, by the Tits lemma,
every two reduced words can be obtained form each other by a sequence of elementary
$2$- and $3$-moves (see \cite[Section~2.1]{ssvz});
by \cite[Theorem~3.5]{ssvz}, every such move either leaves
$(F(\ii), \tilde B(\ii))$ unchanged, or replaces it by an adjacent seed.
The non-simply-laced case is not yet fully understood.
\end{remark}

\begin{remark}
\label{rem:finite-generation-double-cells}
In view of Theorem~\ref{double-cell=upper-bound-explicit}, the
upper cluster algebra $\overline\AA(\ii)$ associated with the
matrix $\tilde B (\ii)$ is always finitely generated,
for it is isomorphic to a $\ZZ$-form of the coordinate ring
of a quasi-affine algebraic variety~$G^{u,v}$.
It is easy to see that not all of these algebras have an acyclic seed;
thus, acyclicity is \emph{not necessary} for finite generation of an upper
cluster
algebra.
\end{remark}

\begin{remark}
\label{rem:tp-bases}
As shown in \cite[Theorem~1.11]{fz-double}, the collection of
minors $F(\ii)$ provides a \emph{total positivity criterion}
in~$G^{u,v}$, in the following sense: an element $x \in G^{u,v}$ is
totally nonnegative if and only if $\Delta(x)>0$ for every
$\Delta \in F(\ii)$.
Theorem~\ref{double-cell=upper-bound-explicit} generates 
lots of other total positivity criteria, obtained from $F(\ii)$ by
replacing the cluster of all $\ii$-exchangeable minors with an
arbitrary cluster in $\AA(\ii)$.
\end{remark}

Combining Theorem~\ref{double-cell=upper-bound-explicit} with
Theorem~\ref{th:acyclic-sufficient-upper}, we obtain a sufficient
condition under which the coordinate ring $\CC[G^{u,v}]$ can be
identified with the cluster algebra $\AA(\ii)_\CC$ defined
(over the complex numbers) by the same
seed $\Sigma(\ii)=(\xx,\tilde B(\ii))$ of geometric type
that was used to define the upper cluster algebra~$\overline\AA(\ii)_\CC$.

\begin{corollary}
\label{col:Guv=cluster-algebra-acyclic}
Let $\ii$ be a reduced word for $(u,v)\in W\times W$
such that the
matrix $B(\ii)$ is mutation equivalent to an acyclic matrix.
Then the isomorphism
\hbox{$\varphi: \FFcal_\CC \to \CC(G^{u,v})$} in
Theorem~\ref{double-cell=upper-bound-explicit} restricts to an
isomorphism of algebras \hbox{${\AA(\ii)}_\CC \to \CC[G^{u,v}]$.}
\end{corollary}

\begin{example}
\label{example:SL3-D4}
In our running example of the reduced word $\ii=(1,2,1,-1,-2,-1)$
for $u=v=w_\circ\in S_3$,
the principal submatrix~$B(\ii)$ is given by
\[
B(\ii)=
\left[
\begin{array}{rrrr}
  0 & -1 &  1 &  0\\
  1 &  0 & -1 &  1\\
 -1 &  1 &  0 & -1\\
  0 & -1 &  1 &  0
\end{array}
\right].
\]
Applying the matrix mutation $\mu_2$ (or~$\mu_3$) results in an
acyclic matrix;
the corresponding transformation of the graph~$\Gamma(B)$ is shown in
Figure~\ref{fig:A2-D4}.
Thus, Corollary~\ref{col:Guv=cluster-algebra-acyclic} applies,
and the coordinate ring of the Bruhat cell $G^{w_\circ,w_\circ}\subset
SL_3(\CC)$ is naturally endowed with a cluster algebra structure.
Furthermore, since the diagram on the right in Figure~\ref{fig:A2-D4} is
an orientation of a Dynkin diagram of type~$D_4$,
it follows that the cluster algebra in question has (finite)
type~$D_4$.
By the results in~\cite{fz-clust2},
such an algebra has $16$~cluster variables, which form $50$~clusters,
each consisting of $4$~variables.
In this particular example, the $16$~cluster variables are:
\begin{itemize}
\item[(i)]
$14$ (of the $19$ available) minors of a $3\times 3$ matrix~$x=(x_{ij})$, namely,
all excluding $\det(x)$ and the four minors
in~\eqref{eq:4-frozen-minors};
\item[(ii)]
the two regular functions
$x_{12}x_{21}x_{33}-x_{12}x_{23}x_{31}-x_{13}x_{21}x_{32}+x_{13}x_{22}x_{31}$
and
$x_{11}x_{23}x_{32}-x_{12}x_{23}x_{31}-x_{13}x_{21}x_{32}+x_{13}x_{22}x_{31}$.
\end{itemize}
We note that each of the $50$ clusters of this cluster algebra
provides a total positivity criterion
(cf.~Remark~\ref{rem:tp-bases}):
a $3\times 3$ matrix~$x$ (not necessarily of determinant~$1$)
is totally positive (that is, all its minors are
positive) if and only if the following $9$~functions take positive
values at~$x$:
\begin{itemize}
\item
the four cluster variables in an arbitrary fixed cluster;
\item
the four minors in~\eqref{eq:4-frozen-minors};
\item
the determinant~$\det(x)$. 
\end{itemize}
The $34$ of these $50$ criteria are shown in
\cite[Figure~8]{fz-double}.
These are precisely the criteria that do not involve the two
``exotic'' cluster variables in~(ii):
each of those $34$ criteria tests $9$ minors of a $3\times 3$ matrix.
\end{example}

\begin{figure}[ht]
\setlength{\unitlength}{1.5pt}
\begin{picture}(60,22)(0,-2)
\put(0,0){\line(1,0){60}}
\put(0,0){\line(1,1){20}}
\put(20,20){\line(1,-1){20}}
\put(20,20){\line(2,-1){40}}

\put(20,0){\vector(1,0){1}}
\put(50,0){\vector(-1,0){1}}
\put(10,10){\vector(-1,-1){1}}
\put(30,10){\vector(-1,1){1}}
\put(40,10){\vector(2,-1){1}}

\put(0,0){\circle*{2}}
\put(40,0){\circle*{2}}
\put(60,0){\circle*{2}}
\put(20,20){\circle*{2}}

\end{picture}
\begin{picture}(30,22)(-15,0)
\put(0,15){\makebox(0,0)
{$\stackrel{\textstyle \mu_2}{\longleftrightarrow}$}}
\end{picture}
\begin{picture}(60,22)(0,-2)
\put(0,0){\line(1,1){20}}
\put(20,20){\line(1,-1){20}}
\put(60,0){\line(-2,1){40}}

\put(10,10){\vector(1,1){1}}
\put(30,10){\vector(1,-1){1}}
\put(40,10){\vector(-2,1){1}}

\put(0,0){\circle*{2}}
\put(40,0){\circle*{2}}
\put(60,0){\circle*{2}}
\put(20,20){\circle*{2}}

\end{picture}

\caption{Transformation of the graph~$\Gamma(B)$  for
  $\CC[G^{w_\circ,w_\circ}]$ in type~$A_2$.}
\label{fig:A2-D4}
\end{figure}

Following \cite[Section~7]{fz-clust2}, an integer square
matrix $B$ is called \emph{$2$-finite} if it is totally
sign-skew-symmetric, and every matrix $B'=(b'_{ij})$
mutation equivalent to $B$ satisfies
$|b'_{k\ell} b'_{k\ell}| \leq 3$ for all $k$ and~$\ell$.
According to \cite[Theorem~1.8]{fz-clust2}, this condition
characterizes cluster algebras of finite type.
In particular, every $2$-finite matrix is mutation equivalent to an
acyclic one, and we obtain the following result.

\begin{corollary}
\label{col:Guv=cluster-algebra-finite}
The conclusion of Corollary~\ref{col:Guv=cluster-algebra-acyclic}
holds whenever the matrix $B(\ii)$ is $2$-finite.
\end{corollary}

\begin{remark}
\label{rem:upper-bound-exact}
It would be interesting to classify
\begin{itemize}
\item
all pairs $(u,v)$ for which
the upper bound $\overline\AA(\ii)_\CC$ in
Theorem~\ref{double-cell=upper-bound-explicit} is equal to the cluster
algebra $\AA(\ii)_\CC$;
\item
all pairs $(u,v)$
for which there exists a reduced word $\ii$ such that $B(\ii)$
is mutation equivalent to an acyclic matrix.
\end{itemize}
By Corollary~\ref{col:Guv=cluster-algebra-acyclic},
the first class of pairs contains the second one.
We suspect that these two classes are different, i.e., there exist
double cells for which $\overline\AA(\ii) = \AA(\ii)$
although this cluster algebra is not acyclic.
In fact, it is not hard to see that in type~$A_r$, the upper bound is
exact, i.e., $\overline\AA(\ii) = \AA(\ii)$,
for each of the cells $G^{e, w_\circ}$ and $G^{w_\circ, w_\circ}$.
(This follows from the arguments in Remark~\ref{rem:ii-dependence}
and the fact that all fundamental representations of $G = SL_{r+1}$
are minuscule.)
Our preliminary calculations indicate that in type~$A_5$, the cluster algebra
associated with $G^{e, w_\circ}$ does \emph{not} have an acyclic seed.
\end{remark}

\subsection{Double Bruhat cells associated with Coxeter elements}

\begin{example}
\label{ex:double-rank-0}
Suppose that $|\Supp(u,v)| = \ell(u) + \ell(v)$, i.e.,
for some (equivalently, any)
reduced word $\ii$ for $(u,v)$, the indices
$|i_1|, \dots, |i_{\ell(u) + \ell(v)}|$ are distinct .
Then the double cell $G^{u,v}$ is a complex torus
of dimension $m = r + \ell (u) + \ell (v)$, so its coordinate ring
is the algebra of Laurent polynomials in the minors of the
family~$F(\ii)$, which in this case consists of the minors
$\Delta_{\omega_j, v^{-1} \omega_j}$ for $j \in [1,r]$, and
$\Delta_{u \omega_i, \omega_i}$ for $i \in \Supp(u,v)$.
\end{example}

\begin{example}
\label{ex:double-Gcc}
Let $c = s_1 \cdots s_r$ be a Coxeter element in~$W$.
(Since our numbering of simple roots is arbitrary, any Coxeter
element could be obtained this~way.)
Consider the double Bruhat cell~$G^{c,c}$.
Pick the following reduced word for~$(c,c)$:
\[
\ii= (-1, \dots, -r, 1, \dots, r).
\]
Applying Definition~\ref{def:edges}, we see that the graph $\Gamma(B(\ii))$ consists
of $r$ disconnected vertices, i.e., $B(\ii)$ is the zero $r \times r$ matrix.
In the terminology of~\cite{fz-clust2}, the corresponding cluster
algebra $\AA$ is of the finite type $A_1^r=A_1\times\cdots\times A_1$ ($r$~times).
The family $F(\ii)$ consists of the $r$ cluster
variables $\Delta_{c \omega_j, c^{-1} \omega_j}$ for $j \in
[1,r]$, together with the $2r$ ``coefficients"
$\Delta_{\omega_i, c^{-1} \omega_i}$ and $\Delta_{c \omega_i, \omega_i}$
for $i \in [1,r]$.
In view of Definition~\ref{def:tildeB-double-cells}, the exchange
relations from the initial cluster take the following form:
\begin{equation}
\label{eq:exchange-Gcc}
\Delta_{c \omega_j, c^{-1} \omega_j} \Delta'_j =
\Delta_{\omega_j, c^{-1} \omega_j} \Delta_{c \omega_j, \omega_j} +
\prod_{i > j} \Delta_{\omega_i, c^{-1} \omega_i}^{-a_{ij}} \cdot
\prod_{i < j} \Delta_{c\omega_i, \omega_i}^{-a_{ij}}.
\end{equation}
Comparing (\ref{eq:exchange-Gcc}) with the
identity~\cite[(1.25)]{fz-double}
(in which we need to replace $i$ by~$j$, and set $u = s_1 \cdots s_{j-1}$
and $v = s_r \cdots s_{j+1}$), we conclude that
$\Delta'_j = \Delta_{\omega_j, \omega_j}$.
Combining Corollary~\ref{col:Guv=cluster-algebra-acyclic}
with Corollaries~\ref{cor:finite-generation-sufficient}
and \ref{cor:cluster-algebra-finite-generation},
we obtain the following result.

\begin{proposition}
For a Coxeter element $c\in W$,
the coordinate ring $\CC[G^{c,c}]$ is a cluster algebra of
type~$A_1^r$.
It has a total of $2r$ cluster variables $\Delta_{c \omega_j, c^{-1} \omega_j}$
and $\Delta_{\omega_j, \omega_j}$, for $j \in [1,r]$, and is
generated by them together with the functions
$\Delta_{\omega_i, c^{-1} \omega_i}^{\pm 1}$ and
$\Delta_{c \omega_i, \omega_i}^{\pm 1}$, for $i \in [1,r]$.
The ideal of relations is generated by the relations~\eqref{eq:exchange-Gcc}.
\end{proposition}

The double cell $G^{c,c}$ appeared in \cite{reshet} where it was
studied from the point of view of integrable systems.
It would be interesting to relate the results in \cite{reshet}
with the cluster algebra structure described above.
\end{example}

\begin{example}
\label{ex:double-Gccinverse}
Let $c$ be a Coxeter element as in Example~\ref{ex:double-Gcc}.
We now consider the double cell $G^{c,c^{-1}}$.
As a reduced word for $(c,c^{-1})$, we choose
$\ii= (-1, \dots, -r, r, \dots, 1)$.
The matrix $B(\ii) = (b_{ij})_{i,j \in [1,r]}$ is then given as follows:
\begin{equation}
\label{eq:Coxeter-bij}
b_{ij} =
\begin{cases}
-a_{ij} & \text{if $i < j$;} \\[.05in]
a_{ij} & \text{if $i > j$.}
\end{cases}
\end{equation}
Thus, the graph $\Gamma(B(\ii))$ is an orientation of the Dynkin graph of~$G$.
In the terminology of~\cite{fz-clust2}, the corresponding cluster
algebra $\AA$ is of finite type, and its type is the
Cartan-Killing type of $G$.
The family $F(\ii)$ consists of the $r$ cluster
variables $\Delta_{c \omega_j, c \omega_j}$ for $j \in
[1,r]$, together with the $2r$ ``coefficients"
$\Delta_{\omega_i, c \omega_i}$ and $\Delta_{c \omega_i, \omega_i}$
for $i \in [1,r]$.
The exchange
relations from the initial cluster take the following form:
\begin{equation}
\label{eq:exchange-Gccinverse}
\Delta_{c \omega_j, c \omega_j} \Delta'_j =
\Delta_{\omega_j, c \omega_j} \Delta_{c \omega_j, \omega_j}
\prod_{i > j} \Delta_{c \omega_i, c \omega_i}^{-a_{ij}} +
\prod_{i > j} (\Delta_{\omega_i, c \omega_i} \Delta_{c\omega_i, \omega_i})^{-a_{ij}} \cdot
\prod_{i < j} \Delta_{c\omega_i, c \omega_i}^{-a_{ij}}.
\end{equation}
In contrast with Example~\ref{ex:double-Gcc}, the
functions $\Delta'_j \in \CC[G^{c,c^{-1}}]$ are no longer minors
in general.
Another dissimilarity with Example~\ref{ex:double-Gcc} is that the
total supply of cluster variables is now much larger.
On the other hand, the same argument as above shows that $\CC[G^{c,c^{-1}}]$
is generated by the $2r$ cluster variables $\Delta_{c \omega_j, c \omega_j}$
and $\Delta'_j$ for $j \in [1,r]$, together with the functions
$\Delta_{\omega_i, c\omega_i}^{\pm 1}$ and
$\Delta_{c \omega_i, \omega_i}^{\pm 1}$ for $i \in [1,r]$.
The ideal of relations among them is generated by the
relations~(\ref{eq:exchange-Gccinverse}).

If 
$\ell (c^2) = 2r$, then the same considerations
show that each of $\CC[G^{c^2,e}]$ and $\CC[G^{e,c^2}]$ is a cluster
algebra of finite type, and its type is the Cartan-Killing type of~$G$.
\end{example}

\begin{remark}
The above results about cluster algebra structures on double
Bruhat cells are likely to extend to Kac-Moody groups.
For instance, Example~\ref{ex:double-Gccinverse} provides a concrete
geometric realization
of all cluster algebras of finite type
(for a special choice of coefficients);
we expect this construction to extend to Kac-Moody groups,
thus providing a geometric realization for a cluster
algebra with an arbitrary acyclic exchange matrix.
\end{remark}

\subsection{Double Bruhat cells associated with base affine spaces}
\label{sec:base-affine}

Let us now consider the pair $(u,v) = (e,w_\circ)$.
The corresponding double cell $G^{e,w_\circ}$ is naturally identified
with the open subset of the base affine space $N_-\backslash G$
given by the conditions $\Delta_{\omega_i, \omega_i} \neq 0$ and
$\Delta_{\omega_i, w_\circ \omega_i} \neq 0$ for all $i$.
Assume that the Cartan matrix $A$ of $G$ is indecomposable,
i.e., the Dynkin graph is connected.
The case where $G = SL_2$ is of type $A_1$ is a special case
of Example~\ref{ex:double-rank-0}: $G^{e,w_\circ}$ is a $2$-dimensional
complex torus, so $\CC[G^{e,w_\circ}]$ is the algebra of Laurent polynomials
in $2$ variables.

\begin{proposition}
\label{pr:GmodN-finite-type}
The list of all simple groups $G$ of rank at least $2$ such that
the cluster algebra associated with $\CC[G^{e,w_\circ}]$ is of finite
type, is given in Table~\ref{tab:GmodN-finite-type}.
\end{proposition}

\begin{table}[ht]
\begin{center}
\begin{tabular}{|c|c|c|c|c|}
\hline
&&&&\\[-.1in]
Cartan-Killing type of $G$ & $A_2$ & $A_3$ & $A_4$ & $B_2$ \\
\hline
&&&&\\[-.1in]
Cluster type of $\CC[G^{e,w_\circ}]$ & $A_1$
& $A_3$
& $D_6$
& $B_2$
 \\[.05in]
\hline
\end{tabular}
\end{center}
\medskip
\caption{Base affine spaces of finite cluster type}
\label{tab:GmodN-finite-type}
\end{table}

\begin{proof}
We apply the results in \cite[Sections~7--9]{fz-clust2}
to the matrix $B(\ii_\circ)$, where $\ii_\circ$ is a reduced word
for $w_\circ$.
Thus, we only have to identify the types of $G$ for which the
corresponding matrix $B(\ii_\circ)$ is $2$-finite.

In view of \cite[Section~8]{fz-clust2}, one can determine whether a
given matrix $B=(b_{k\ell})$ is $2$-finite by examining its
\emph{diagram}, defined as the directed graph~$\Gamma(B)$
in which each edge $(k,\ell)$ is assigned the weight~$|b_{k\ell}
b_{\ell k}|$. The matrix mutations descend onto the level of such
diagrams, and a diagram
corresponds to a $2$-finite matrix if and only if all diagrams in its
mutation equivalence class have edge weights~$\leq 3$.
(See \cite[Section~8]{fz-clust2} for further details.)
In the rest of Section~\ref{sec:double-cells}, we find it more
convenient to use the language of diagram
mutations instead of matrix mutations.
The reader may choose to ignore this switch, and stick with the
matrix setup, which would cause no harm.


We use the following choice for $\ii_\circ$
(see, e.g., \cite[Exercise~V.\S6.2]{bourbaki}).
Since we assume that $A$ is indecomposable,
the Dynkin graph is of one of the Cartan-Killing types
$A_r, B_r, \dots, G_2$ (see, e.g., \cite{bourbaki}).
In particular, it is a tree, hence a bipartite graph.
Thus, $[1,r]$ is the disjoint union of the two subsets
$I_+$ and $I_-$ none of which contains a pair of indices $i$ and $j$
such that $a_{ij} < 0$.
Let $h$ denote the Coxeter number of $\Phi$, i.e., the order of
any Coxeter element in~$W$.
Let $\ii_+$ (resp., $\ii_-$) be some permutation of indices from
$I_+$ (resp., $I_-$).
Then the word
\begin{equation}
\label{eq:special-rw}
\ii_\circ\stackrel{\rm def}{=}
\underbrace{\ii_- \ii_+ \ii_- \cdots \ii_\mp
 \ii_\pm}_h
\end{equation}
(concatenation of $h$ segments)
is a reduced word for~$w_\circ$.
Recall that $h$ is even for all
types except $A_r$ with $r$ even;
in the exceptional case of type $A_{2e}$, we have $h = 2e +1$.

First we show that the diagram $\Gamma(B(\ii_\circ))$
is $2$-finite in each of the cases in
Table~\ref{tab:GmodN-finite-type}.
Note that the number of vertices of
$\Gamma(B(\ii_\circ))$ is $n = \ell (w_\circ) - r$.
For $G$ of type $A_2$, we have $n = 3 -2 = 1$, so the
corresponding cluster type is $A_1$.
For $G$ of type $A_3$ or $B_2$, we have
$w_\circ = c^2$ for a Coxeter element $c$, and $\ell (w_\circ) = 2r$.
Thus, the corresponding cluster type is the same as the
Cartan-Killing type of $G$, by the last remark in
Example~\ref{ex:double-Gccinverse}.
For $G$ of type $A_4$, we have $\ii_\circ = (1,3,2,4,1,3,2,4,1,3)$
(in the standard labeling of simple roots), so $n = 10 - 4 = 6$.
The corresponding diagram $\Gamma(B(\ii_\circ))$ is shown on the
left in Figure~\ref{fig:A4-D6}
(with all edge weights equal to~$1$).
Applying the mutations $\mu_1$ and $\mu_2$ at the two leftmost
vertices, we obtain the diagram on the right in
Figure~\ref{fig:A4-D6}.
Since this diagram is an orientation of the Dynkin diagram $D_6$,
the case $A_4$ is done.

\begin{figure}[ht]
\setlength{\unitlength}{1.5pt}
\begin{picture}(100,65)(0,-2)
\put(0,0){\line(1,0){80}}
\put(0,0){\line(2,1){20}}
\put(40,20){\line(-1,1){20}}
\put(20,40){\line(2,1){20}}
\put(40,20){\line(3,1){20}}
\put(20,40){\line(1,0){80}}
\put(40,20){\line(2,-1){20}}
\put(60,60){\line(2,-1){20}}

\put(0,0){\vector(1,0){72}}
\put(40,20){\vector(-2,-1){32}}
\put(40,20){\vector(-1,1){15}}
\put(60,60){\vector(-2,-1){32}}
\put(100,40){\vector(-3,-1){52}}
\put(20,40){\vector(1,0){72}}
\put(80,0){\vector(-2,1){32}}
\put(100,40){\vector(-2,1){32}}

\put(0,0){\circle*{2}}
\put(80,0){\circle*{2}}
\put(40,20){\circle*{2}}
\put(20,40){\circle*{2}}
\put(100,40){\circle*{2}}
\put(60,60){\circle*{2}}

\end{picture}
\begin{picture}(50,65)(-25,0)
\put(0,35){\makebox(0,0)
{$\stackrel{\textstyle\mu_1 \mu_2}{\longleftarrow\!\!\longrightarrow}$}}
\end{picture}
\begin{picture}(100,65)(0,-2)
\put(0,0){\line(1,0){80}}
\put(40,20){\line(-2,-1){20}}
\put(40,20){\line(-1,1){20}}
\put(60,60){\line(-2,-1){20}}
\put(20,40){\line(1,0){80}}

\put(80,0){\vector(-1,0){72}}
\put(0,00){\vector(2,1){32}}
\put(20,40){\vector(1,-1){15}}
\put(20,40){\vector(2,1){32}}
\put(100,40){\vector(-1,0){72}}

\put(0,0){\circle*{2}}
\put(80,0){\circle*{2}}
\put(40,20){\circle*{2}}
\put(20,40){\circle*{2}}
\put(100,40){\circle*{2}}
\put(60,60){\circle*{2}}

\end{picture}
\caption{Diagrams for $\CC[G^{e,w_\circ}]$ in type $A_4$.}
\label{fig:A4-D6}
\end{figure}

To complete the argument, it suffices to show that, in each of the
types not appearing in Table~\ref{tab:GmodN-finite-type}, the
diagram $\Gamma(B(\ii_\circ))$ contains a subdiagram which is an
``extended Dynkin tree diagram'' in the sense of
\cite[Proposition~9.3]{fz-clust2}.
This is done by a direct inspection.
For example, in type $A_5$,
the reduced word $\ii_\circ$ is given by
$\ii_\circ = (1,3,5,2,4,1,3,5,2,4,1,3,5,2,4)$,
and the diagram $\Gamma(B(\ii_\circ))$
is shown in Figure~\ref{fig:A5-E7hat};
removing vertices $4$ and $5$, we obtain a subdiagram
of affine type $E_7^{(1)}$ which is an extended Dynkin tree diagram
(its vertices are enlarged in
Figure~\ref{fig:A5-E7hat}).

\begin{figure}[ht]
\setlength{\unitlength}{1.5pt}
\begin{picture}(180,85)(0,-2)
\put(0,0){\line(1,0){100}}
\put(60,20){\line(1,0){100}}
\put(20,40){\line(1,0){100}}
\put(80,60){\line(1,0){100}}
\put(40,80){\line(1,0){100}}
\put(0,0){\line(3,1){20}}
\put(100,0){\line(3,1){20}}
\put(20,40){\line(2,-1){20}}
\put(120,40){\line(2,-1){20}}
\put(20,40){\line(3,1){20}}
\put(120,40){\line(3,1){20}}
\put(40,80){\line(2,-1){20}}
\put(140,80){\line(2,-1){20}}
\put(60,20){\line(2,-1){20}}
\put(60,20){\line(3,1){20}}
\put(80,60){\line(2,-1){20}}
\put(80,60){\line(3,1){20}}

\put(0,0){\vector(1,0){92}}
\put(60,20){\vector(1,0){92}}
\put(20,40){\vector(1,0){92}}
\put(80,60){\vector(1,0){92}}
\put(40,80){\vector(1,0){92}}

\put(60,20){\vector(-3,-1){48}}
\put(160,20){\vector(-3,-1){48}}
\put(60,20){\vector(-2,1){32}}
\put(160,20){\vector(-2,1){32}}
\put(80,60){\vector(-3,-1){48}}
\put(180,60){\vector(-3,-1){48}}
\put(80,60){\vector(-2,1){32}}
\put(180,60){\vector(-2,1){32}}
\put(100,0){\vector(-2,1){32}}
\put(120,40){\vector(-3,-1){48}}
\put(120,40){\vector(-2,1){32}}
\put(140,80){\vector(-3,-1){48}}

\put(0,0){\circle*{4}}
\put(100,0){\circle*{4}}
\put(60,20){\circle*{2}}
\put(160,20){\circle*{4}}
\put(20,40){\circle*{4}}
\put(120,40){\circle*{4}}
\put(80,60){\circle*{2}}
\put(180,60){\circle*{4}}
\put(40,80){\circle*{4}}
\put(140,80){\circle*{4}}

\end{picture}

\caption{A diagram for $\CC[G^{e,w_\circ}]$ in type $A_5$.}
\label{fig:A5-E7hat}
\end{figure}

Similarly, in type $D_4$, we have 
$\ii_\circ = (1,3,4,2,1,3,4,2,1,3,4,2)$ (with $2$ being the branching point
of $D_4$), and the diagram $\Gamma(B(\ii_\circ))$ has vertices
$1,\dots, 8$ corresponding to the first $8$ terms in
$\ii_\circ$; removing vertices $5$, $6$ and $7$, we obtain a subdiagram
of affine type $D_4^{(1)}$ which is an extended Dynkin tree diagram.
The straightforward check in all the other cases is left to the
reader.
\end{proof}

\begin{remark}
\label{rem:leclerc-connection}
As shown by B.~Leclerc~\cite{leclerc}, the types in
Table~\ref{tab:GmodN-finite-type} are precisely those satisfying
the ``Berenstein--Zelevinsky conjecture'' on the multiplicative
structure of the dual canonical basis.
This is by no means a coincidence; in fact, one of the
authors of the present paper (A.Z.) suggested that
counterexamples to this conjecture must exist in all cases where the cluster
algebra associated with $\CC[G^{e,w_\circ}]$ is of infinite type.
The relationships between cluster algebras and dual canonical
bases will be explored in a separate publication.
\end{remark}

\pagebreak[2]

\subsection{Open double Bruhat cells}
\label{sec:open-double-bruhat}

Our last example is the pair $(u,v) = (w_\circ,w_\circ)$.
By Proposition~\ref{pr:double-cell-by-eqs},
the corresponding double cell $G^{w_\circ,w_\circ}$ is
the open subset of $G$
given by the non-vanishing conditions $\Delta_{w_\circ \omega_i, \omega_i} \neq 0$ and
$\Delta_{\omega_i, w_\circ \omega_i} \neq 0$ for all~$i\in[1,r]$.
Again, we assume that the Cartan matrix $A$ of $G$ is indecomposable,

\begin{proposition}
\label{pr:G-finite-type}
If $G$ is of type $A_1$ (resp.,~$A_2$) then
the cluster algebra associated with $\CC[G^{w_\circ,w_\circ}]$ is of finite
type $A_1$ (resp.,~$D_4$).
For every other simple group $G$, this cluster algebra is of infinite type.
\end{proposition}

\begin{figure}[ht]
\setlength{\unitlength}{1.5pt}
\begin{picture}(160,42)(0,-2)
\put(0,0){\line(1,0){120}}
\put(40,20){\line(1,0){120}}
\put(20,40){\line(1,0){120}}
\put(80,40){\line(-2,-1){20}}
\put(40,20){\line(-2,-1){20}}
\put(20,40){\line(1,-1){40}}
\put(100,20){\line(2,1){20}}
\put(100,20){\line(1,-1){20}}
\put(120,0){\line(2,1){20}}
\put(140,40){\line(1,-1){20}}

\put(60,0){\vector(1,0){52}}
\put(100,20){\vector(1,0){52}}
\put(80,40){\vector(1,0){52}}
\put(60,0){\vector(-1,0){52}}
\put(100,20){\vector(-1,0){52}}
\put(20,40){\vector(1,-1){15}}
\put(40,20){\vector(1,-1){15}}
\put(160,20){\vector(-1,1){15}}
\put(120,0){\vector(-1,1){15}}
\put(0,0){\vector(2,1){32}}
\put(40,20){\vector(2,1){32}}
\put(140,40){\vector(-2,-1){32}}
\put(160,20){\vector(-2,-1){32}}
\put(80,40){\vector(-1,0){52}}

\put(0,0){\circle*{4}}
\put(60,0){\circle*{2}}
\put(120,0){\circle*{4}}
\put(40,20){\circle*{4}}
\put(100,20){\circle*{4}}
\put(160,20){\circle*{2}}
\put(20,40){\circle*{4}}
\put(80,40){\circle*{2}}
\put(140,40){\circle*{4}}
\end{picture}

\caption{A diagram for $\CC[G^{w_\circ,w_\circ}]$ in type $A_3$.}
\label{fig:A3-D5hat}
\end{figure}
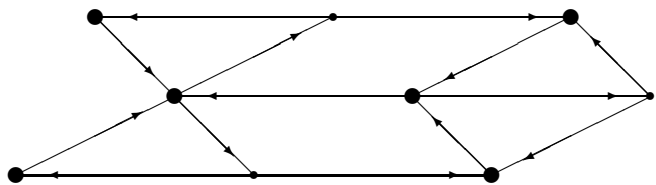

\begin{proof}
The proof is similar to that of
Proposition~\ref{pr:GmodN-finite-type}.
For type $A_1$, the claim is a special case of each of the
Examples~\ref{ex:double-Gcc} and \ref{ex:double-Gccinverse}.
The type $A_2$ has already been treated in
Example~\ref{example:SL3-D4}.

As before, to finish the proof it suffices to show that, in each
of the remaining Cartan-Killing types, the diagram
$\Gamma(B(\ii))$ contains an extended Dynkin tree subdiagram.
This is again done by direct inspection.
For example, in type $A_3$, the diagram for the reduced word
$\ii = (-1,-3,-2,-1,-3,-2,1,3,2,1,3,2)$ is shown in
Figure~\ref{fig:A3-D5hat}.
Removing from this diagram the vertices
$4$, $5$ and $9$ counting from the left, we obtain a subdiagram of
type $D_5^{(1)}$ (its vertices are enlarged in
Figure~\ref{fig:A3-D5hat}).
\end{proof}

\section{Proof of Proposition~\ref{pr:coprime}
}
\label{sec:proof-prop-coprime}


\begin{lemma}
\label{lem:coprime-criterion}
A seed $\Sigma = (\xx, \tilde B)$ of geometric type is coprime if and only if
no two columns of $\tilde B$ are proportional to each other
with the proportionality coefficient being a ratio of two odd integers.
\end{lemma}

\begin{proof}
According to Definition~\ref{def:coprime-seed}, $\Sigma$~is coprime
if and only if
any common divisor of $P_i$ and $P_j$ ($i\neq j$) belongs to the
set~$\PP$ of Laurent monomials in $x_{n+1},\dots,x_m$.
Since both $P_i$ and $P_j$ have no monomial factors,
this condition is equivalent to $P_i$ and $P_j$ being coprime as
polynomials in $x_1, \dots, x_m$ with integer coefficients.
Thus, for the rest of the proof, we are treating $P_1,\dots,P_n$ as
elements of $\ZZ[x_1,\dots,x_m]$.

Let $\tilde B_j$
denote the $j$th
column of~$\tilde B$.
It suffices to show that $P_j$ and $P_k$ have a 
common factor if and only if
$\tilde B_k=\pm\frac{b}{a}\tilde B_j\,$,
where $a$ and $b$ are odd and coprime positive integers.
The ``if'' part is clear: if $\tilde B_k=\pm\frac{b}{a}\tilde B_j\,$,
with $a$ and $b$ as stated,
then $P_j=L^a+M^a$ and $P_k=L^b+M^b$, for some monomials $L$ and~$M$
 (see~(\ref{eq:exchange})),
and $L+M$ is a common factor of $P_j$ and~$P_k\,$.

Let us prove the ``only if'' part.
For any $P\in\ZZ[x_1,\dots,x_m]$, we denote by $N(P)\subset\RR^m$ the
Newton polytope of~$P$, i.e., the convex hull in $\RR^m$ of all lattice points
$(a_1, \dots, a_m)$ such that the coefficient of the monomial
$x_1^{a_1} \cdots x_m^{a_m}$ in $P$ is nonzero.
We note that the Newton polytope
$N(PQ)$ of a product is equal to
the Minkowski sum $N(P) + N(Q)$.

Since $P_j$ is a sum of two monomials,
$N(P_j)$ is a line segment;
it is furthermore parallel to~$\tilde B_j\,$.
It follows that, for every non-monomial factor $P$ of $P_j$,
the polytope $N(P)$ is also a line segment parallel to~$N(P_j)$.

Suppose that $P_j$ and $P_k$ are not coprime.
Then $N(P_j)$ and $N(P_k)$ are collinear,
and so
$\tilde B_k=\pm\frac{b}{a}\tilde B_j$,
with $a$ and $b$ coprime positive integers.
Thus, $P_j=L^a+M^a$ and $P_k=L^b+M^b$, for some monomials $L$ and~$M$
in disjoint sets of variables.
It remains to show that if $L^a+M^a$ and $L^b+M^b$
have a non-trivial common factor
then both $a$ and $b$ are odd.

First of all, we notice that the polynomials
$t^a+1$ and $t^b + 1$ have a common complex root if and only if
$a$ and $b$ are odd (since the roots of $t^a + 1$ are of the form
$t = e^{mi\pi/a}$ for some odd integer $m$).
Therefore, if one of $a$ and $b$ is even then
$t^a + 1$ and $t^b + 1$ are coprime elements of the polynomial ring
$\CC[t]$.
It follows that there exist polynomials $f(t),g(t)\in \CC[t]$ such that
$f(t)(1+t^a)+g(t)(1+t^b)=1$.
Substituting $t=M/L$ and clearing denominators, we obtain
$F \cdot (L^a+M^a) + G \cdot (L^b+M^b)=L^c$, for some
$F,G\in\CC[x_1,\dots,x_m]$ and
some integer~$c\geq 0$.
It follows that
a common factor of $L^a+M^a$ and~$L^b+M^b$
must be a monomial in $x_1, \dots, x_m$, and we are done.
\end{proof}

Proposition~\ref{pr:coprime} is an immediate consequence of
Lemma~\ref{lem:coprime-criterion} together with the following
statement, which previously appeared in~\cite[Lemma~1.2]{gsv}.

\begin{lemma}
\label{lem:rank-mutation-invariant}
Matrix mutations preserve the rank of $\tilde B$.
\end{lemma}

\begin{proof}
Fix an index $k \in [1,n]$ and a sign $\varepsilon \in \{\pm 1\}$.
The rule (\ref{eq:mutation})
describing a matrix mutation in direction~$k$
can be rewritten as follows:
\begin{equation}
\label{eq:matrix-mutation-restated}
b'_{ij} =
\begin{cases}
-b_{ij} & \text{if $i=k$ or $j=k$;} \\[.05in]
b_{ij} + b_{ik} \max(0, \varepsilon b_{kj})
+ \max(0, - \varepsilon b_{ik}) b_{kj} & \text{otherwise.}
\end{cases}
\end{equation}
To see this, substitute
\[
\max(0, \varepsilon b_{kj}) =
\frac{\varepsilon b_{kj} + |b_{kj}|}{2}, \quad
\max(0, -\varepsilon b_{ik}) =
\frac{-\varepsilon b_{ik} + |b_{ik}|}{2} \, .
\]
To prove the lemma, observe that (\ref{eq:matrix-mutation-restated})
can be restated as
\begin{equation}
\label{eq:mutation-product}
\tilde B' = (J_{m,k} + E_k) \,\tilde B\, (J_{n,k} + F_k)
=J_{m,k}\,\tilde B J_{n,k} + J_{m,k}\,\tilde B F_k + E_k \,\tilde B J_{n,k} \, ,
\end{equation}
where
\begin{itemize}
\item
$J_{m,k}$ (resp., $J_{n,k}$) denotes the diagonal $m\times m$ (resp.,~$n\times
n$) matrix whose diagonal entries are all~$1$'s, except for $-1$ in
the $k$th position;
\item
$E_k$
is the $m\times m$ matrix whose only nonzero entries are
$e_{ik}\!=\!\max(0, -\varepsilon b_{ik})$;
\item
$F_k$
is the $n\times n$ matrix whose only nonzero entries are
$f_{kj}\!=\!\max(0, \varepsilon b_{kj})$.
\end{itemize}
(Here we used that $b_{ii}=0$ for all~$i$,
and therefore $E_k\,\tilde B F_k=0$.)
Since both $J_{m,k}+E_k$ and $J_{n,k}+F_k$ are invertible, it follows that
$\tilde B'$ and $\tilde B$ have the same rank.
Lemma~\ref{lem:rank-mutation-invariant} and
Proposition~\ref{pr:coprime} are proved.
\end{proof}

\section{Proof of Theorem~\ref{th:upper-bound-universal}}
\label{sec:proof-upper-bound}

We retain the notation of Section~\ref{sec:main-results}.
In particular, $P_1,\dots,P_n$ are the exchange polynomials
(see~(\ref{eq:exchange})) for a seed $\Sigma=(\xx,\pp,B)$.

Our first two lemmas hold for any seed
(not necessarily coprime).

\begin{lemma}
\label{lem:upper-bound-intersection-1}
The upper bound at an arbitrary seed can be expressed as follows:
\begin{equation}
\label{eq:upper-bound-2}
\Upper (\Sigma) = \bigcap_{j = 1}^n
\ZP [x_1^{\pm 1}, \ldots, x_{j-1}^{\pm 1}, x_j, x'_j,
 x_{j+1}^{\pm 1}, \ldots, x_{n}^{\pm 1}].
\end{equation}
\end{lemma}

\begin{proof}
It is enough to show that
\begin{equation}
\label{eq:two-laurent-conditions-2}
\ZP [\xx^{\pm 1}] \cap \ZP [\xx_1^{\pm 1}] =
\ZP [x_1,x'_1, x_2^{\pm 1}, \dots, x_n^{\pm 1}].
\end{equation}
The ``$\supseteq$" inclusion is clear.
To prove~``$\subseteq$", note that
any $y\in\ZP [\xx^{\pm 1}]$ is of the form
$$y = \sum_{m = -N}^N c_m x_1^m,$$
where $c_m\in\ZP [x_2^{\pm 1}, \dots, x_n^{\pm 1}]$ for all~$m$.
Substituting $x_1 = P_1/x'_1$, we get
$$y = \sum_{m = 0}^N c_m P_1^m {x'_1}^{-m} +
\sum_{m = 1}^N (c_{-m}/P_1^m) \cdot {x'_1}^{m}.$$
If, in addition, $y \in \ZP [\xx_1^{\pm 1}]$, then
$c_{-m}/P_1^m \in \ZP [x_2^{\pm 1}, \dots, x_n^{\pm 1}]$
for $m = 1, \dots, N$.
Hence
$$y = \sum_{m = 0}^N c_m x_1^m + \sum_{m = 1}^N (c_{-m}/P_1^m) \cdot {x'_1}^{m}
\in \ZP [x_1,x'_1, x_2^{\pm 1}, \dots, x_n^{\pm 1}],$$
as desired.
\end{proof}

Note for future reference that the above proof of
(\ref{eq:two-laurent-conditions-2}) implies the following
description of the subalgebra
$\ZP [x_1,x'_1, x_2^{\pm 1}, \dots, x_n^{\pm 1}]$
inside the ambient field $\FFcal$.

\begin{lemma}
\label{lem:two-laurent-conditions-1}
An element $y \in \FFcal$ belongs to
$\ZP [x_1,x'_1, x_2^{\pm 1}, \dots, x_n^{\pm 1}]$ if and only
if $y$ can be written in the form $y = \sum_{m = -N}^N c_m x_1^m$,
where each coefficient $c_m$ belongs to
$\ZP [x_2^{\pm 1}, \dots, x_n^{\pm 1}]$, and
$c_{-m}$ is divisible by $P_1^m$ in $\ZP [x_2^{\pm 1}, \dots, x_n^{\pm 1}]$
for $m = 1, \dots, N$.
\end{lemma}

To prove Theorem~\ref{th:upper-bound-universal}, we will need the
following refinement of (\ref{eq:upper-bound-2}), which requires
some coprimality assumptions.

\begin{proposition}
\label{pr:upper-bound-intersection-2}
Suppose that $n \geq 2$, and $P_1$ is coprime with $P_j$ for
$j = 2, \dots, n$.
Then
\begin{equation}
\label{eq:upper-bound-intersection-2}
\Upper (\Sigma) = \bigcap_{j = 2}^n
\ZP [x_1, x'_1, x_2^{\pm 1}, \ldots, x_{j-1}^{\pm 1}, x_j, x'_j,
 x_{j+1}^{\pm 1}, \ldots, x_{n}^{\pm 1}].
\end{equation}
\end{proposition}

\begin{proof}
Comparing (\ref{eq:upper-bound-intersection-2}) with
(\ref{eq:upper-bound-2}), we see that it is enough to show
that the term with $j = 2$ in 
(\ref{eq:upper-bound-intersection-2}) is equal to the intersection
of the two terms with $j = 1$ and $j = 2$ in 
(\ref{eq:upper-bound-2}).
Let us ``freeze" the cluster variables $x_3, \dots, x_n$, i.e., let us
view $\ZP [ x_3^{\pm 1}, \dots, x_n^{\pm 1}]$ as the new ground
ring (cf.~\cite[Proposition~2.6]{fz-clust1}).
In other words, we can assume without loss of generality that
$n=2$ from the start, so that
(\ref{eq:upper-bound-intersection-2})
reduces to
\begin{equation}
\label{eq:rk-2-upper=lower}
\ZP [x_1,x'_1, x_2^{\pm 1}] \cap \ZP [x_1^{\pm 1}, x_2,x'_2] =
\ZP [x_1,x'_1, x_2,x'_2].
\end{equation}

The proof of (\ref{eq:rk-2-upper=lower}) breaks into two cases.

\noindent \emph{Case 1:} $b_{12} = b_{21} = 0$.
The exchange relations (\ref{eq:exchange}) take the form
$$x_1 x'_1 = P_1 \in\ZP, \qquad
x_2 x'_2 = P_2\in\ZP \, .$$
Then the algebra
\[
\ZP [x_1,x'_1, x_2,x'_2]
=\ZP[x_1,x_2]+\ZP[x_1',x_2]+\ZP[x_1,x_2']+\ZP[x_1',x_2']
\]
consists of all Laurent polynomials
\[
y = \sum_{m_1, m_2 \in \ZZ} c_{m_1,m_2} x_1^{m_1} x_2^{m_2}
\]
with coefficients in $\ZP$ such that
\begin{align}
& \text{$c_{m_1,m_2}$ is divisible by $P_1^{- m_1}$ whenever
$m_1 < 0 \leq m_2$;} \\
& \text{$c_{m_1,m_2}$ is divisible by $P_2^{- m_2}$ whenever
$m_2 < 0 \leq m_1$;} \\
\label{eq:divisible-by-P1P2}
&\text{$c_{m_1,m_2}$ is divisible by $P_1^{- m_1}P_2^{- m_2}$ whenever
$m_1 < 0$ and $m_2 < 0$.}
\end{align}
On the other hand, using Lemma~\ref{lem:two-laurent-conditions-1},
we see that 
$\ZP [x_1,x'_1, x_2^{\pm 1}] \cap \ZP [x_1^{\pm 1}, x_2,x'_2]$
has a similar
description: the only difference is that (\ref{eq:divisible-by-P1P2})
gets replaced by the condition
that $c_{m_1,m_2}$ is divisible by $P_1^{- m_1}$ and by $P_2^{- m_2}$ whenever
$m_1 < 0$ and $m_2 < 0$.
The equivalence of this condition to (\ref{eq:divisible-by-P1P2}) is
ensured by the assumption that $P_1$ and $P_2$ are coprime.

\smallskip

\noindent \emph{Case 2:} $b_{12}b_{21} < 0$.
Now the exchange relations (\ref{eq:exchange}) take the form
\begin{eqnarray}
\label{eq:x1-x'1}
&&x_1 x'_1 = P_1 = q_2 x_2^c + r_2,\\
\label{eq:x2-x'2}
&&x_2 x'_2 = P_2 = q_1 x_1^b + r_1,
\end{eqnarray}
for some $q_1, q_2, r_1, r_2 \in\PP$, where both $b = |b_{12}|$ and $c= |b_{21}|$
are positive integers.

\begin{lemma}
\label{lem:rk-2-upper=lower-1}
$\ZP [x_1, x_2^{\pm 1}] \cap \ZP [x_1^{\pm 1}, x_2,x'_2] =
\ZP [x_1,x_2,x'_2]$.
\end{lemma}

\begin{proof}
The ``$\supseteq$" inclusion is obvious.
To prove ``$\subseteq$'', assume that
$y\in\ZP [x_1, x_2^{\pm 1}] \cap \ZP [x_1^{\pm 1},x_2,x'_2]$.
Applying if necessary the relation (\ref{eq:x2-x'2}), we can write $y$
in the form
\begin{equation}
\label{eq:x1-expansion-rank2-cm}
y = \sum_{m\in\ZZ} x_1^m (c_m + c_m'(x_2) + c_m''(x_2'))\,,
\end{equation}
where $c_m\!\in\!\ZP$, and $c_m'$ and $c_m''$ are polynomials over~$\ZP$
without constant term.
(All but finitely many terms in (\ref{eq:x1-expansion-rank2-cm})
vanish.)
Let $M$ be the smallest integer~such that the coefficient of $x_1^M$
in (\ref{eq:x1-expansion-rank2-cm}) is nonzero.
If $M\!\geq\! 0$, then we are done.
Otherwise, expanding
\begin{equation*}
\textstyle y = \sum 
x_1^m (c_m + c_m'(x_2) + c_m''(\textstyle\frac{q_1 x_1^b + r_1}{x_2}))
\end{equation*}
as a Laurent polynomial in $x_1$ and~$x_2$, we see that the terms with the smallest power of~$x_1$
in this expansion are
$x_1^M (c_M + c_M'(x_2) + c_M''(\textstyle\frac{r_1}{x_2}))\neq 0$,
contradicting the condition $y\in\ZP [x_1, x_2^{\pm 1}]$.
\end{proof}

\begin{lemma}
\label{lem:Z-in-Z+Z}
$
\ZP[x_1,x_1',x_2^{\pm 1}] =
\ZP[x_1,x_1',x_2,x_2'] + \ZP[x_1,x_2^{\pm 1}]
$.
\end{lemma}

\begin{proof}
The ``$\supseteq$" inclusion is obvious.
The only difficulty in proving ``$\subseteq$'' is to show that
${x'_1}^N x_2^{-M} \in \ZP[x_1,x_1',x_2,x_2'] + \ZP[x_1,x_2^{\pm 1}]$
for all $M,N > 0$.
Let $p = - q_1/r_1$.
Then (\ref{eq:x2-x'2}) can be rewritten as
$$x_2^{-1} - p x_1^b x_2^{-1} = r_1^{-1}x_2'\, ,$$
implying
$$x_2^{-1} \equiv p x_1^b x_2^{-1} \equiv p^2 x_1^{2b} x_2^{-1}
\equiv \cdots \equiv p^N x_1^{Nb} x_2^{-1} \bmod
\ZP[x_1,x_2'] \, .$$
Now write $x_2^{-1} = p^N x_1^{Nb} x_2^{-1} + y$
with $y \in \ZP[x_1,x_2']$, and raise this equality to the
$M^{\rm th}$ power to obtain
$$x_2^{-M} \in \ZP[x_1,x_2'] + x_1^N \ZP[x_1,x_2^{-1}] \,.$$
Hence
${x'_1}^N x_2^{-M} \in \ZP[x_1,x'_1,x_2'] + \ZP[x_1,x_2^{\pm 1}]$,
as desired.
\end{proof}

Now everything is ready for the proof of (\ref{eq:rk-2-upper=lower})
in Case~2.
Using Lemmas~\ref{lem:rk-2-upper=lower-1} and~\ref{lem:Z-in-Z+Z}
together with the inclusion
$\ZP[x_1,x_1',x_2,x_2'] \subset \ZP[x_1^{\pm 1},x_2,x_2']$,
we obtain:
\begin{align*}
\,&\ZP[x_1,x_1',x_2^{\pm 1}]\,\cap\, \ZP[x_1^{\pm 1},x_2,x_2']\\
=\,
&(\ZP[x_1,x_1',x_2,x_2'] + \ZP[x_1,x_2^{\pm 1}])
\cap \ZP[x_1^{\pm 1},x_2,x_2']\\
=\,
&\ZP[x_1,x_1',x_2,x_2'] + (\ZP[x_1,x_2^{\pm 1}]
\cap \ZP[x_1^{\pm 1},x_2,x_2'])\\
=\,
&\ZP[x_1,x_1',x_2,x_2'] + \ZP [x_1,x_2,x'_2]\\
=\,
&\ZP[x_1,x_1',x_2,x_2'].
\end{align*}
Equality (\ref{eq:rk-2-upper=lower}) and
Proposition~\ref{pr:upper-bound-intersection-2} are proved.
\end{proof}

Turning to the proof of Theorem~\ref{th:upper-bound-universal},
we can assume without loss of generality that $\Sigma' = \Sigma_1$
is obtained from the seed $\Sigma$ by mutation in direction $1$.
If $n=1$, there is nothing to prove, so we assume that $n \geq 2$.
Comparing (\ref{eq:upper-bound-intersection-2}) with its counterpart for
$\Upper (\Sigma_1)$, we see that it suffices to show the
following.

\begin{lemma}
\label{lem:reduction-to-rank-2}
Let $x_2''$ be the cluster variable exchanged with $x_2$ under the
corresponding mutation of~$\Sigma_1$.
Then
\begin{equation}
\label{eq:reduction-to-rank-2}
\ZP [x_1,x'_1, x_2,x'_2, x_3^{\pm 1}, \dots, x_n^{\pm 1}]
= \ZP [x_1,x'_1, x_2,x''_2, x_3^{\pm 1}, \dots, x_n^{\pm 1}].
\end{equation}
\end{lemma}

\begin{proof}
As in the proof of
Proposition~\ref{pr:upper-bound-intersection-2},
we can ``freeze" the
cluster variables $x_3, \dots, x_n$, thus
reducing the proof of (\ref{eq:reduction-to-rank-2})
to the case~$n=2$.
In that case, (\ref{eq:reduction-to-rank-2}) takes the form
\begin{equation}
\label{eq:rk-2-upper-bound}
\ZP [x_1,x'_1, x_2,x'_2]
= \ZP [x_1,x'_1, x_2, x''_2].
\end{equation}

By symmetry, 
it is enough to show that
$x''_2 \in \ZP [x_1,x'_1, x_2,x'_2]$.
Consider the same two cases as in the proof of
Proposition~\ref{pr:upper-bound-intersection-2}.
In Case~1, we have
$x''_2 = px'_2$ for some $p \in \PP$, which makes
the inclusion $x''_2 \in \ZP [x_1,x'_1, x_2,x'_2]$ obvious.
In Case~2, the exchange relations from the initial seed take
the form (\ref{eq:x1-x'1})--(\ref{eq:x2-x'2}), and the variable $x''_2$
is obtained by
\begin{equation}
\label{eq:x2-x''2}
x_2 x''_2 = q_3 {x'_1}^b + r_3;
\end{equation}
using the mutation rule (\ref{eq:p-mutation}), we obtain
the relation
\begin{equation}
\label{eq:rk2-coeff-relation}
r_1 r_3 = q_1 q_3 r_2^b
\end{equation}
(cf.\ \cite[(2.12)]{fz-clust1}).
It follows that
\begin{align*}
x''_2
&=x_2^{-1}(q_3 (x_1')^b+r_3)\\
&=x_2^{-1} q_3 (x_1')^b r_1^{-1}(x_2 x'_2-q_1 x_1^b)+r_3 x_2^{-1}\\
&=q_3 r_1^{-1} (x_1')^b x_2'
  - x_2^{-1}(q_1 q_3 r_1^{-1}(x_1 x_1')^b -r_3)\\
&=q_3 r_1^{-1} (x_1')^b x_2'
  - q_1 q_3 r_1^{-1} ((q_2 x_2^c+r_2)^b - r_2^b)/x_2\\
&\in \ZP [x_1,x'_1, x_2,x'_2].
\end{align*}
This completes the proof of (\ref{eq:rk-2-upper-bound}),
Lemma~\ref{lem:reduction-to-rank-2}, and
Theorem~\ref{th:upper-bound-universal}.
\end{proof}

\begin{remark}
\label{rem:rank2-Thm1.18}
The equality (\ref{eq:rk-2-upper=lower}), which plays the central
role in the above argument, is a special case of
Theorem~\ref{th:acyclic-sufficient-upper}
for the cluster algebras of rank~$2$;
the acyclicity of~$\Sigma$ is automatic in rank~$2$.
\end{remark}

\section{Proof of
 Theorem~\ref{th:acyclic-sufficient-lower}}
\label{sec:proof-acyclic-sufficient}

We start with the proof of the ``if" part: if a seed $\Sigma$
is acyclic then the standard monomials are linearly
independent over~$\ZP$.
We note that Definition~\ref{def:acyclic-seed} can be restated as
follows: a sign-skew-symmetric matrix $B$ is acyclic if and only if there exists
a permutation $\sigma \in S_n$ such that
$b_{\sigma (i), \sigma (j)} \geq 0$ for $i > j$.
(This follows at once from a well-known fact that every partial
order admits a linear extension.)
Renumbering if necessary the elements of the initial cluster
$\xx$, we can thus assume that $b_{ij} \geq 0$ for~$i>j$.

Let us label the standard monomials in $x_1,x_1',\dots,x_n,x_n'$
by the lattice points
$\mathbf{m} = (m_1, \dots, m_n) \in \ZZ^{n}$ by setting
\begin{equation}
\label{eq:standard-monomial-notation-1}
\xx^{\langle\mathbf{m}\rangle} = x_1^{\langle m_1\rangle} \cdots
x_n^{\langle m_n\rangle},
\end{equation}
where we abbreviate
\begin{equation}
\label{eq:standard-monomial-notation-2}
x_j^{\langle m_j\rangle}=
\begin{cases}
x_j^{m_j} & \text{if $m_j \geq 0$;} \\[.05in]
(x'_j)^{-m_j} & \text{if $m_j < 0$.}
\end{cases}
\end{equation}
(Note that each $x_j^{\langle m_j\rangle}$, and
hence~$\xx^{\langle\mathbf{m}\rangle}$,  is a Laurent polynomial
in $x_1,\dots,x_n\,$.)
We order both the standard monomials
$\xx^{\langle\mathbf{m}\rangle}$ and the ordinary Laurent
monomials $\xx^{\mathbf{m}} = x_1^{m_1} \cdots x_n^{m_n}$
lexicographically as follows:
\begin{equation}
\label{eq:order-monomial}
\text{$\mathbf{m} \prec \mathbf{m'}$ if the first
nonzero difference $m'_j - m_j$ is positive.}
\end{equation}
For example,
\[
x_1^{-1}\prec x_2^{-1}\prec \cdots\prec x_n^{-1}\prec 1
\prec x_n\cdots\prec x_2\prec x_1\,.
\]
Since $b_{ij} \geq 0$ for $i > j$, the (lexicographically)
first Laurent monomial that appears in
$x_j^{\langle m_j\rangle}$ is equal to
$x_j^{m_j}$ times some monomial in $x_{j+1}, \dots, x_n$.
This implies the following:
\begin{equation}
\label{eq:leading-monomial-1}
\text{if $\mathbf{m} \prec \mathbf{m}'$, then the first monomial
in $\xx^{\langle \mathbf{m} \rangle}$ precedes the one
in $\xx^{\langle \mathbf{m'} \rangle}$.}
\end{equation}
The linear independence of standard monomials over
$\ZP$ follows at once from~(\ref{eq:leading-monomial-1}).

We will deduce the ``only if" part of
Theorem~\ref{th:acyclic-sufficient-lower} from the following
proposition.

\begin{proposition}
\label{pr:acyclic-necessary-2}
Suppose a seed $\Sigma = (\xx, \pp, B)$ is such that
$1 \to 2 \to \cdots \to \ell \to 1$ is an oriented cycle
in~$\Gamma (B)$.
Then
\begin{equation}
\label{eq:x1-xl=sum}
x'_1 \cdots x'_\ell
=\sum_{K\subsetneq\{1,\dots,\ell\}} f_K(x_1,\dots,x_n) \prod_{k\in K}
x_k'\,,
\end{equation}
for some polynomials $f_K(x_1,\dots,x_n)\in\ZP[\xx] = \ZP[x_1, \dots, x_n]$.
\end{proposition}

The proof of Proposition~\ref{pr:acyclic-necessary-2}
utilizes the following algebraic identity.

\begin{lemma}
Let $I$ be a finite non-empty set, let $i \mapsto i^+$ be a cyclic
permutation of~$I$, and let $(u_i)_{i\in I}$ and $(v_i)_{i\in I}$
be two families of (commuting) indeterminates.
Then
\begin{equation}
\label{eq:diffcomb}
\doublesubscript{\sum}{J\subset I}{J\cap J^+=\emptyset} (-1)^{|J|}
\prod_{i\in I- (J\cup J^+)} (u_i+v_i) \cdot\prod_{j\in J} (u_j
v_{j^+})
=
\prod_{i \in I} u_i + \prod_{i \in I} v_i \ .
\end{equation}
\end{lemma}

\begin{proof}
Expanding the left-hand side of (\ref{eq:diffcomb}),
we rewrite it as follows:
$$\sum_{K \subset I} c_K \,\prod_{i \in K} u_i \,\prod_{j \in I - K} v_j,$$
where the coefficients $c_K$ are given by
$$c_K = \sum_{J \subset \{i \in K:\, i^+ \in I-K\}} (-1)^{|J|} \ .$$
We see that $c_K = 0$ unless the set $\{i \in K: i^+ \in I-K\}$ is empty.
Since $i \mapsto i^+$ is a cyclic permutation of $I$, the only two choices for $K$
with $c_K \neq 0$ are $K = \emptyset$ and $K = I$.
In both these case, $c_K = 1$, and (\ref{eq:diffcomb}) follows.
\end{proof}

\noindent
\textbf{Proof of Proposition~\ref{pr:acyclic-necessary-2}.}
We apply (\ref{eq:diffcomb}) in the following situation:
take $I = \{1, \dots, \ell\}$; set $j^+ = j+1$ for $j < \ell$, and
$\ell^+ = 1$; finally, let $u_j = M_j^-/x_j$ and $v_j=M_j^+/x_j$,
where $M_j^+$ and $M_j^-$ are the two terms on the right-hand side
of the exchange relation~(\ref{eq:exchange}).
Under these settings, we have $u_j + v_j = x'_j$, and so the
term on the left in (\ref{eq:diffcomb}) that corresponds to
$J = \emptyset$ turns into $x'_1 \cdots x'_\ell$.
It remains to show that all products $u_j v_{j^+}$ as well as
$u_1 \cdots u_\ell$ and $v_1 \cdots v_\ell$ are regular monomials
in $x_1, \dots, x_n$ (i.e., they do not contain negative powers).
For $u_j v_{j^+}$, it suffices to
note that $b_{j^+,j} < 0$ for all $j = 1, \dots, \ell$, so
the Laurent monomial $u_j$ (resp., $v_{j^+}$) contains a positive power of
$x_{j^+}$ (resp., $x_j$).
Similarly, the denominator $x_{j^+}$ (resp., $x_j$) of the Laurent
monomial $u_{j^+}$ (resp.,~$v_j$) is absorbed by the term
$u_j$ (resp.,~$v_{j^+}$) in the product $u_1 \cdots u_\ell$
(resp.,~$v_1 \cdots v_\ell$).
\endproof

To complete the proof of Theorem~\ref{th:acyclic-sufficient-lower},
we need to show that, under the assumptions of
Proposition~\ref{pr:acyclic-necessary-2},
the standard monomials in $x_1,x_1',\dots,x_n,x_n'$ are linearly
dependent over~$\ZP$.
This follows from Proposition~\ref{pr:acyclic-necessary-2} by a
straightening argument: repeatedly applying the exchange
relations~\eqref{eq:exchange} to the right-hand side
of~\eqref{eq:x1-xl=sum}, we can transform it into a linear combination
of standard monomials none of which contains $x'_1 \cdots x'_\ell$.
\endproof

\section{Proof of
 Theorem~\ref{th:acyclic-sufficient-upper}}
\label{sec:proof-acyclic-sufficient-upper}

The proof will require some preparation.
For the first several statements,
we only assume that a seed $\Sigma = (\xx, \pp, B)$ is acyclic;
the coprimality condition will be invoked later.
As in Section~\ref{sec:proof-acyclic-sufficient}, we assume
without loss of generality that $b_{ij} \geq 0$ for $i>j$.

Let $\ZP^{\rm st} [x_2,x'_2,\ldots, x_n, x'_n]$
(resp.,~$\ZP^{\rm st} [x_1,x_2,x'_2,\ldots, x_n, x'_n]$)
denote the $\ZZ\PP$-linear span (resp.,~$\ZZ\PP[x_1]$-linear span) of 
the standard monomials in $x_2, x'_2, \ldots, x_n, x'_n$.
As above, by using (\ref{eq:exchange}) as straightening
relations, we conclude that
\begin{equation}
\label{eq:standard-inclusion}
\ZP [x_2, x'_2, \ldots, x_n, x'_n]
\subset \ZP^{\rm st} [x_1,x_2,x'_2,\ldots, x_n, x'_n] .
\end{equation}

\begin{definition}
Let
$\varphi:\ZP [x_2,x'_2,\ldots, x_n, x'_n]
\to \ZP [x_{2}^{\pm 1}, \ldots, x_{n}^{\pm 1}]$
denote the algebra homomorphism obtained as a composition
$$\ZP [x_2,x'_2,\ldots, x_n, x'_n] \to
\ZP [x_1, x_{2}^{\pm 1}, \ldots, x_{n}^{\pm 1}]
\to \ZP [x_{2}^{\pm 1}, \ldots, x_{n}^{\pm 1}],$$
where the first map is the embedding given by
(\ref{eq:exchange}), and the second map is the specialization
$x_1 \mapsto 0$.
\end{definition}

\begin{lemma}
\label{lem:phi-1}
There is a direct sum decomposition
$$\ZP [x_2,x'_2,\ldots, x_n, x'_n] = {\rm Ker} (\varphi)
\oplus \ZP^{\rm st} [x_2,x'_2,\ldots, x_n, x'_n] \ .$$
\end{lemma}

\begin{proof}
The decomposition
\[
\ZP [x_2,x'_2,\ldots, x_n, x'_n] = {\rm Ker} (\varphi)
+ \ZP^{\rm st} [x_2,x'_2,\ldots, x_n, x'_n]
\]
follows from (\ref{eq:standard-inclusion}).
To prove that the sum is direct, we need to show that
the restriction of $\varphi$ to
$\ZP^{\rm st} [x_2,x'_2,\ldots, x_n, x'_n]$
is injective.
The proof is similar to that of the ``if" part
of Theorem~\ref{th:acyclic-sufficient-lower} given in
Section~\ref{sec:proof-acyclic-sufficient}.
We label the standard monomials in $x_2,x'_2,\ldots, x_n, x'_n$
by the lattice points $\mathbf{m} = (m_2, \dots, m_n) \in \ZZ^{n-1}$
using the conventions
(\ref{eq:standard-monomial-notation-1})--(\ref{eq:standard-monomial-notation-2}),
and order the standard monomials
$\xx^{\langle\mathbf{m}\rangle}$ and the ordinary Laurent
monomials $\xx^{\mathbf{m}}$ according to (\ref{eq:order-monomial}).
As in Section~\ref{sec:proof-acyclic-sufficient}, we use the condition that
$b_{ij} \geq 0$ for $i > j$ to conclude that the (lexicographically)
first Laurent monomial appearing in
$\varphi (x_j^{\langle m_j\rangle})$ is
$x_j^{m_j}$ times some monomial in $x_{j+1}, \dots, x_n$.
This implies that
\begin{equation}
\label{eq:leading-monomial}
\text{if $\mathbf{m} \prec \mathbf{m}'$, then the first monomial
in $\varphi (\xx^{\langle \mathbf{m} \rangle})$ precedes the one
in $\varphi (\xx^{\langle \mathbf{m'} \rangle})$.}
\end{equation}
The injectivity of $\varphi$ on the span of standard monomials
is immediate from~(\ref{eq:leading-monomial}).
\end{proof}

\begin{definition}
\label{def:LT1}
For a Laurent polynomial $y\in\ZP [x_{1}^{\pm 1}, \ldots, x_{n}^{\pm
 1}]$, let $\operatorname{LT}_1(y)$ denote the leading term of~$y$
with respect to~$x_1$.
More precisely, $\operatorname{LT}_1(y)$ is the
 sum of all Laurent monomials involving the \emph{smallest} power of
 $x_1$ which appear in the expansion of~$y$ with nonzero coefficient.
\end{definition}

\begin{lemma}
\label{lem:leading-term-x1}
Suppose that
\begin{equation}
\label{eq:x1-expansion-standard}
y = \sum_{m = a}^b c_m x_1^m, \quad c_m \in
\ZP^{\rm st} [x_2, x'_2, \ldots, x_n, x'_n] \,, \quad c_a\neq 0.
\end{equation}
Then $\operatorname{LT}_1(y)=\varphi(c_a) x_1^a$.
\end{lemma}

\begin{proof}
It suffices to check that $\varphi(c_a)\neq 0$.
This follows from the injectivity of $\varphi$ on $\ZP^{\rm st} [x_2,
x'_2, \ldots, x_n, x'_n]$, which holds by Lemma~\ref{lem:phi-1}.
\end{proof}

We now generalize Lemma~\ref{lem:rk-2-upper=lower-1}.

\begin{lemma}
\label{lem:upper=lower-1}
We have
\begin{equation}
\label{eq:upper-bound-inductive-2}
\ZP [x_1, x_{2}^{\pm 1}, \ldots, x_{n}^{\pm 1}]
\cap \ZP [x_{1}^{\pm 1}, x_2,x'_2,\ldots, x_n, x'_n]
= \ZP [x_1, x_2, x'_2, \ldots, x_n, x'_n].
\end{equation}
\end{lemma}

\begin{proof}
The ``$\supseteq$" inclusion is obvious.
To prove the ``$\subseteq$" part, consider an element
$y\neq 0$ of the left-hand side of (\ref{eq:upper-bound-inductive-2}).
Since $y \in \ZP [x_{1}^{\pm 1}, x_2,x'_2,\ldots, x_n, x'_n]$,
the inclusion (\ref{eq:standard-inclusion})
implies that $y$ can be expressed in the
form~(\ref{eq:x1-expansion-standard}).
By Lemma~\ref{lem:leading-term-x1},
$\operatorname{LT}_1(y)=\varphi(c_a) x_1^a$.
Now $y\in\ZP [x_1, x_{2}^{\pm 1}, \ldots, x_{n}^{\pm 1}]$
implies that $a\geq 0$, and (\ref{eq:x1-expansion-standard})
shows that $y\in\ZP [x_1, x_2, x'_2, \ldots, x_n, x'_n]$, as desired.
\end{proof}

The next step is the hardest part of the proof.
We obtain the following description of the image of~$\varphi$.

\begin{lemma}
\label{lem:phi-image}
$\operatorname{Im}(\varphi)=\ZP [x_2, x_{2}^{(-)}, \ldots, x_n, x_{n}^{(-)}]$,
where we use the notation
\begin{equation}
\label{eq:x-(minus)}
x_j^{(-)}=
\begin{cases}
x'_j & \text{if $b_{1j}=0$;} \\[.05in]
x_j^{-1} & \text{if $b_{1j}\neq 0$.}
\end{cases}
\end{equation}
\end{lemma}

\begin{proof}
For $j = 2, \dots, n$, we set
$$M_j = x_j^{-1} \prod_{i > j} x_i^{b_{ij}}.$$
(Recall that, by the acyclicity condition, the exponents $b_{ij}$ on
the right-hand side are nonnegative.)
Let $J = \{j \in [2,n]: b_{1j} = 0\}$.
In this notation, we have
$$\varphi (x'_j)=
\begin{cases}
x'_j & \text{if $j \in J$;} \\[.05in]
M_j & \text{if $j \in [2,n]-J$.}
\end{cases}
$$
Comparing this with the definition of $x_j^{(-)}$
in~(\ref{eq:x-(minus)}), 
we conclude that
$${\rm Im} (\varphi) \subseteq
\ZP [x_2, x_{2}^{(-)}, \ldots, x_n, x_{n}^{(-)}].$$
To prove the reverse inclusion, we only need to show that
$x_j^{-1} \in {\rm Im} (\varphi)$ for $j \in [2,n]-J$.
This will require some preparation.

Every Laurent monomial $\xx^{\mathbf{m}}\in\ZP[x_2^{\pm
 1},\dots,x_n^{\pm 1}]$
can be uniquely written as
$$\xx^{\mathbf{m}} = \mathbf{M}^{\mathbf{\ell}} =
\prod_{2\leq j\leq n} M_j^{\ell_j},
$$
where the tuples of exponents $\mathbf{m} = (m_2, \dots, m_n)$
and $\mathbf{\ell} = (\ell_2, \dots, \ell_n)$ are related by
\begin{equation}
\label{eq:l-to-m}
m_j = - \ell_j + \sum_{2\leq i < j} b_{ji} \ell_i \, .
\end{equation}
(Again, the acyclicity condition ensures that the coefficients
$b_{ji}$ on the right-hand side are nonnegative.)
Now consider the multiplicative monoid
\[
\text{$\mathcal{M} = \{\xx^{\mathbf{m}} = \mathbf{M}^{\mathbf{\ell}}:
\ell_k \geq 0$ for $k \in [2,n]$, $m_j \geq 0$ for $j \in J\}$.}
\]
In view of (\ref{eq:l-to-m}), we have
$\mathbf{M}^{\mathbf{\ell}} \in \mathcal{M}$ if and only if
\begin{equation}
\label{eq:M-conditions}
\ell_k \geq 0 \quad (k \in [2,n]), \quad
\ell_j \leq \sum_{2\leq i < j} b_{ji} \ell_i \quad
(j \in J).
\end{equation}
It follows readily from (\ref{eq:l-to-m})
and (\ref{eq:M-conditions}) that
$x_j^{-1} \in \mathcal{M}$ for $j \in [2,n]-J$.
(The inequalities $\ell_k\geq 0$ are proved by induction on~$k$.)
Thus, to prove the desired assertion that
$x_j^{-1} \in {\rm Im} (\varphi)$ for $j \in [2,n]-J$,
it is enough to show that
$\mathcal{M} \subset {\rm Im} (\varphi)$.

Let $M = \mathbf{M}^{\mathbf{\ell}} \in \mathcal{M}$.
We prove that $M \in {\rm Im} (\varphi)$ by induction
on ${\rm deg}(M) = \sum_{k=2}^n \ell_k$.
The base $M = 1$ is trivial, so let us assume
that ${\rm deg}(M) > 0$, and that all elements of
$\mathcal{M}$ whose degree is less than $\deg(M)$ belong to~${\rm Im} (\varphi)$.

Let $j$ be the maximal index such that $\ell_j >0$.
As an easy consequence of (\ref{eq:M-conditions}), we have $M/M_j \in \mathcal{M}$.
Hence $M/M_j \in {\rm Im} (\varphi)$ by the induction assumption.
If $j \in [2,n]-J$, then
$M_j = \varphi (x'_j) \in {\rm Im} (\varphi)$, and
$M = (M/M_j) \cdot M_j \in {\rm Im} (\varphi)$, as desired.
So suppose that $j \in J$.
Let us denote
\[
M_j^- = x_j^{-1} \prod_{2\leq k < j} x_k^{|b_{kj}|}.
\]
In this notation, the exchange relation
(\ref{eq:exchange}) can be written as
$x'_j = p_j^+ M_j + p_j^- M_j^-\,.$
(Here we use that $x_1$ does not contribute to the right-hand side
because~$j\in J$.)
Multiplying both sides by $M/M_j$, we get
\[
(M/M_j) \cdot x'_j = p_j^+ M + p_j^- M \cdot (M_j^-/M_j).
\]
Since $(M/M_j) \cdot x'_j \in {\rm Im} (\varphi)$,
to finish the proof it suffices to show that
$M \cdot (M_j^-/M_j) \in {\rm Im} (\varphi)$.
(Here we rely on working over $\ZP$ rather than a subring thereof.)

Since $\ell_j > 0$, it follows from (\ref{eq:M-conditions})
that there exists $i$ between $2$ and $j-1$ such that
$b_{ji} \ell_i >0$.
Fix such an index $i$ and define the monomial
$M' = \prod_{2\leq k\leq n} M_k^{\ell'_k}$ as follows:
\begin{align}
\nonumber
&\text{$\ell'_k = 0$ for $2 \leq k <i$;}
\\
\label{eq:M'}
&\text{$\ell'_i = 1$;}
\\
\nonumber &\text{$\ell'_k = \min(\ell_k, \sum_{2\leq h < k} b_{kh} \ell'_h)$
for $i < k \leq n$.}
\end{align}
(Note that $\ell_k'\geq 0$ for all~$k$.)
Writing $M \cdot (M_j^-/M_j)=(M/M') \cdot (M'M_j^-/M_j)$, we see that
the claim
$M \cdot (M_j^-/M_j) \in {\rm Im} (\varphi)$
is a consequence of the
following two statements:
\begin{eqnarray}
\label{eq:M-over-M'}
&&M/M' \in \mathcal{M};\\
\label{eq:M'MM}
&&\text{$M'M_j^-/M_j = \xx^{\mathbf{m}}$ with
all $m_k \geq 0$ for $k \in [2,n]$.}
\end{eqnarray}
Indeed, ${\rm deg}(M/M') < {\rm deg}(M)$,
(\ref{eq:M-over-M'}), and the inductive assumption imply that
$M/M' \in {\rm Im} (\varphi)$; moreover, (\ref{eq:M'MM})
yields $M'M_j^-/M_j \in \ZP[x_2,\dots,x_n]\subset {\rm Im} (\varphi)$.

To prove (\ref{eq:M-over-M'}), we note that the conditions
(\ref{eq:M-conditions}) for $M/M'$ can be rewritten as
\begin{align}
\label{eq:M-over-M'-conditions-1}
&\ \ \ell'_k \leq \ell_k \quad (k \in [2,n]), \\[.05in]
\label{eq:M-over-M'-conditions-2}
&- \ell'_k+\sum_{2\leq h<k} b_{kh} \ell'_h \leq
- \ell_k+\sum_{2\leq h<k} b_{kh} \ell_h \quad
(k \in J).
\end{align}
The inequalities (\ref{eq:M-over-M'-conditions-1}) are immediate from
(\ref{eq:M'}) and our choice of~$i$.
To prove~(\ref{eq:M-over-M'-conditions-2}), first note that
the right-hand side is nonnegative by~(\ref{eq:M-conditions})
(remember that~$M\in\mathcal{M}$).
Hence (\ref{eq:M-over-M'-conditions-2}) holds for~$k\leq i$.
For $k\geq i$, we have either $\ell_k'=\ell_k$ or $\ell_k'=\sum_{2\leq
 h < k} b_{kh} \ell'_h\,$.
In both cases, (\ref{eq:M-over-M'-conditions-2}) is immediately
checked.

It remains to verify (\ref{eq:M'MM}).
It follows from (\ref{eq:M'}) and (\ref{eq:l-to-m}) that the only
variable that contributes to~$M'$ with negative exponent
(namely,~$-1$) is~$x_i\,$.
The only variable that contributes to~$M_j^-$ with negative exponent
(namely,~$-1$) is~$x_j\,$.
And the only variables that contribute to $M_j^{-1}$ with
(potentially) negative exponents (namely,~$-b_{kj}$)
are the $x_k$ with~$k>j$.
Let us check each of these cases.
First, $m_i = |b_{ij}| - 1 \geq 0$
(note that $b_{ij}<0$ because $b_{ji}>0$ and
$B$ is sign-skew-symmetric).
Second, $m_j=- \ell'_j + \sum_{2\leq h < j} b_{jh} \ell'_h-1+1\geq 0$.
Finally, suppose that~$k > j$.
We note that
$\ell'_k = 0$ because $\ell_k = 0$.
Also, $\ell'_j > 0$ because $\ell_j >0$ and $b_{ji}\ell'_i
= b_{ji} >0$ by our choice of~$i$.
Hence
$m_k = - \ell'_k + \textstyle\sum_{2\leq h < k} b_{kh} \ell'_h - b_{kj}
\geq b_{kj} (\ell'_j -1) \geq 0$,
as desired.
\end{proof}

\noindent
\textbf{Proof of Theorem~\ref{th:acyclic-sufficient-upper}.}
So far, we only relied on the acyclicity of the seed~$\Sigma$.
Now, the condition of $\Sigma$ being coprime comes into play.
We prove Theorem~\ref{th:acyclic-sufficient-upper}
by induction on $n$, the rank of the cluster algebra~$\AA$.
The base case $n=1$ is established in
(\ref{eq:two-laurent-conditions-2}).
So we assume that $n \geq 2$, and that our statement is true for all
algebras of smaller rank.
We note that the rank~2 case was done in
Section~\ref{sec:proof-upper-bound}
(see Remark~\ref{rem:rank2-Thm1.18}).

Using (\ref{eq:upper-bound-2}), we write
\[
\Upper (\Sigma) =
\ZP [x_1,x'_1, x_{2}^{\pm 1}, \ldots, x_{n}^{\pm 1}]
\cap
\bigcap_{j = 2}^n
\ZP [x_1^{\pm 1}, \ldots, x_{j-1}^{\pm 1}, x_j, x'_j,
 x_{j+1}^{\pm 1}, \ldots, x_{n}^{\pm 1}].
\]
Applying the induction assumption to the
seed on the set of indices $[2,n]$ obtained from
$\Sigma$ by ``freezing" the variable $x_1$
(i.e., replacing the ground ring $\ZP$ by $\ZP[x_1^{\pm 1}]$,
cf.\ \cite[Proposition~2.6]{fz-clust1}), we~get:
\[
\bigcap_{j = 2}^n
\ZP [x_1^{\pm 1}, \ldots, x_{j-1}^{\pm 1}, x_j, x'_j,
 x_{j+1}^{\pm 1}, \ldots, x_{n}^{\pm 1}]
=
\ZP [x_{1}^{\pm 1}, x_2,x'_2,\ldots, x_n, x'_n].
\]
The claim $\Lower (\Sigma) = \Upper (\Sigma)$
is thus reduced to proving the equality
\begin{equation}
\label{eq:upper-bound-inductive}
\ZP [x_1,x'_1, x_{2}^{\pm 1}, \ldots, x_{n}^{\pm 1}]
\cap \ZP [x_{1}^{\pm 1}, x_2,x'_2,\ldots, x_n, x'_n]
\!=\! \ZP [x_1, x'_1, \ldots, x_n, x'_n].
\end{equation}

The following lemma plays a crucial role in the proof
of~(\ref{eq:upper-bound-inductive}).

\begin{lemma}
\label{lem:phi-image-saturation}
Suppose that a seed $\Sigma$ is coprime and such that
$b_{ij} \geq 0$ for $i>j$.
Then 
$\operatorname{Im}(\varphi)$ 
has the following
saturation property:
if $z \in \ZP [x_{2}^{\pm 1}, \ldots, x_{n}^{\pm 1}]$
is such that $z P_1 \in \operatorname{Im}(\varphi)$,
then $z \in \operatorname{Im}(\varphi)$.
\end{lemma}

\begin{proof}
Recall the notation $J = \{j \in [2,n]: b_{1j} = 0\}$.
The case $J=\emptyset$ is 
immediate from Lemma~\ref{lem:phi-image},
so we assume that $J\neq\emptyset$.
Now we apply the induction assumption of
Theorem~\ref{th:acyclic-sufficient-upper} to
the seed $\Sigma_1$ on the set of indices $J$
obtained from~$\Sigma$ by forgetting the cluster variable $x_1$, and treating
the variables $x_j$ for $j \in [2,n] - J$ as coefficients.
In view of Lemma~\ref{lem:phi-image}, we have
$$\Lower(\Sigma_1) = \ZP [x_2, x_{2}^{(-)}, \ldots, x_n, x_{n}^{(-)}]
= \operatorname{Im}(\varphi).$$
Using~(\ref{eq:upper-bound-2}), we see that the equality
$\Lower (\Sigma_1) = \Upper (\Sigma_1)$
takes the form
\[
\operatorname{Im}(\varphi)
= \bigcap_{j \in J}
\ZP [x_2^{\pm 1}, \ldots, x_{j-1}^{\pm 1}, x_j, x'_j,
 x_{j+1}^{\pm 1}, \ldots, x_{n}^{\pm 1}].
\]
Therefore, to prove the lemma it suffices to check that each subalgebra
\[
\ZP [x_2^{\pm 1}, \ldots, x_{j-1}^{\pm 1}, x_j, x'_j,
x_{j+1}^{\pm 1}, \ldots, x_{n}^{\pm 1}]
\subset \ZP [x_{2}^{\pm 1}, \ldots, x_{n}^{\pm 1}],
\]
for $j \in J$, satisfies the saturation property.
Let us write a Laurent polynomial
$z \in \ZP [x_{2}^{\pm 1}, \ldots, x_{n}^{\pm 1}]$ in the form
$$z = \textstyle\sum_{m=-N}^N c_m x_j^m, \quad
\quad c_m \in
\ZP [x_2^{\pm 1}, \ldots, x_{j-1}^{\pm 1},
x_{j+1}^{\pm 1}, \ldots, x_{n}^{\pm 1}].$$
By Lemma~\ref{lem:two-laurent-conditions-1},
$z \in \ZP [x_2^{\pm 1}, \ldots, x_{j-1}^{\pm 1}, x_j, x'_j,
x_{j+1}^{\pm 1}, \ldots, x_{n}^{\pm 1}]$ if and only if $c_{-m}$
is divisible by $P_j^m$ for all~$m>0$.
This immediately implies the saturation property in question since
$P_1$ and $P_j$ are coprime.
\end{proof}

We are now ready to prove~(\ref{eq:upper-bound-inductive}).
The ``$\supseteq$" inclusion is obvious.
Let us prove the ``$\subseteq$" part.
Consider an element
$y$ of the left-hand side of~(\ref{eq:upper-bound-inductive}).
If the leading term $\operatorname{LT}_1(y)$
(see Definition~\ref{def:LT1})
involves a nonnegative power of~$x_1$,
then $y\in\ZP [x_1, x_{2}^{\pm 1}, \ldots, x_{n}^{\pm 1}]$,
and we are done by Lemma~\ref{lem:upper=lower-1}.
Otherwise, as in the proof of Lemma~\ref{lem:upper=lower-1},
we express $y$ in the form~(\ref{eq:x1-expansion-standard}), so
by Lemma~\ref{lem:leading-term-x1}, we have
$\operatorname{LT}_1(y)=\varphi(c_a) x_1^a\,$,
where by assumption $a<0$.
Using induction on~$|a|$, to finish the proof
it suffices to find an element
$y' \in \ZP [x_1, x'_1, \ldots, x_n, x'_n]$
such that $\operatorname{LT}_1(y') = \operatorname{LT}_1(y)$.
The inclusion
$y \in \ZP [x_1,x'_1, x_{2}^{\pm 1}, \ldots, x_{n}^{\pm 1}]$
implies by Lemma~\ref{lem:two-laurent-conditions-1}
that $\varphi (c_a)= P_1^{|a|} z$
for some $z \in \ZP [x_{2}^{\pm 1}, \ldots, x_{n}^{\pm 1}]$.
Applying Lemma~\ref{lem:phi-image-saturation} several times, we conclude that
$z = \varphi (c)$ for some $c \in \ZP [x_2, x'_2, \ldots, x_n, x'_n]$.
The element $y' = c \cdot ({x'_1})^{|a|}$ belongs to
$\ZP [x_1, x'_1, \ldots, x_n, x'_n]$, and we have,
by the definition of~$\varphi$,
\[
\operatorname{LT}_1(y')
=\varphi(c)x_1^a P_1^{|a|}
=\varphi (c_a)\, x_1^a
=\operatorname{LT}_1(y),
\]
which completes the proof of Theorem~\ref{th:acyclic-sufficient-upper}.
\qed


\section{Proof of Theorem~\ref{th:lower-bound-exact}}
\label{sec:proof-lower-bound-exact}

To prove the ``if" part of Theorem~\ref{th:lower-bound-exact},
we need to show that acyclicity of
$\Sigma$ implies~$\Lower (\Sigma)=\AA(\Sigma)$,
that is, any cluster variable $x$ of the algebra $\AA = \AA(\Sigma)$ can be written as a
polynomial in $x_1,x_1',\dots,x_n,x_n'$ with coefficients in~$\ZP$.
By Theorem~\ref{th:acyclic-sufficient-upper}, this is true under the
additional assumption that $\Sigma$ is coprime.
In particular, it is true for the cluster algebra
with \emph{universal coefficients} (see \cite[Section~5]{fz-clust1})
defined by the matrix~$B$.
The corresponding formulas expressing cluster variables as polynomials
in $x_1,x_1',\dots,x_n,x_n'$ hold in any cluster algebra,
and the claim follows.

It remains to show that if a seed $\Sigma$ is \emph{not} acyclic
then $\Lower (\Sigma) \neq \AA$.
A more precise statement is given in the following proposition.

\begin{proposition}
\label{pr:acyclic-necessary-1}
Suppose that a seed $\Sigma = (\xx,\pp, B)$ 
is such that
\[
1\to 2 \to \cdots \to \ell \to 1
\]
is an induced oriented cycle
in~$\Gamma (\Sigma)$.
(That is, there are no edges connecting these $\ell$ vertices except for
the $\ell$ edges in the cycle.)
Let $\Sigma^{(0)}, \Sigma^{(1)}, \dots, \Sigma^{(\ell-1)}$
be the sequence of seeds such that $\Sigma^{(0)} = \Sigma$, and
each $\Sigma^{(k)} = (\xx^{(k)}, \pp^{(k)}, B^{(k)})$
is obtained from $\Sigma^{(k-1)}$ by a mutation in
direction $k$ for $k\! = \!1, \dots, \ell \!-\! 1$.
Then the cluster variable $y$ such that $\{y\} = \xx^{(\ell-1)} -
\xx^{(\ell-2)}$ does not belong to $\Lower (\Sigma)$.
\end{proposition}

We prove Proposition~\ref{pr:acyclic-necessary-1} by
constructing a ``valuation'' 
that is nonnegative on
$\Lower (\Sigma)$
but takes a negative value at $y \in \AA$.

To be more specific, recall from Section~\ref{sec:main-results}
that the algebra $\overline\AA$
consists of the elements of $\FFcal$
which are Laurent polynomials with coefficients in $\ZP$
in the cluster elements of any seed which is mutation equivalent
to $\Sigma$.
This algebra contains the cluster algebra $\AA$ (see
Corollary~\ref{cor:three-inclusions}).
Let $\overline\AA_{\rm sf} \subset \overline\AA - \{0\}$
denote the intersection of $\overline\AA$ with the \emph{semifield}
generated by all cluster variables and the coefficient
group~$\PP$.
Thus, $\overline\AA_{\rm sf}$ consists of all nonzero elements of
$\overline\AA$
that can be written as
subtraction-free rational expressions in the cluster variables and
the elements of~$\PP$.
Equivalently, the elements of $\overline\AA_{\rm sf}$ are
those elements of $\overline\AA$ that can be written as ratios of
polynomials in the cluster variables and the elements of $\PP$ with
nonnegative integer coefficients.

\begin{definition}
\label{def:tropical-valuation}
A \emph{tropical valuation} on $\overline\AA$ is a map
$\nu: \overline\AA - \{0\} \to \RR$ that satisfies the following conditions:
\begin{align}
\label{eq:valuation-0}
&\nu(p) = 0 \quad {\rm for} \,\, p \in \ZP; \\
\label{eq:valuation-1}
&\nu(xy) = \nu(x)+\nu(y) \quad {\rm for} \,\, x,y \in \overline\AA - \{0\}; \\
\label{eq:valuation-2}
&\nu(x+y) \geq \min(\nu(x),\nu(y)) \quad {\rm for} \,\, x,y,x+y \in
\overline\AA- \{0\}. \\
\label{eq:valuation-2'}
&\nu(x+y) = \min(\nu(x),\nu(y)) \quad {\rm for} \,\, x,y \in
\overline\AA_{\rm sf}.
\end{align}
\end{definition}

\begin{lemma}
\label{lem:valuation}
For any cluster $\xx = \{x_1, \dots, x_n\}$ of $\AA$ and
any 
$(\nu_1, \dots, \nu_n)\in\RR^n$, there is
a tropical valuation $\nu$ on $\overline\AA$ such that $\nu(x_i) = \nu_i$ for
$i = 1, \dots, n$.
\end{lemma}

\begin{proof}
As in the proof of Lemma~\ref{lem:coprime-criterion},
for a nonzero Laurent polynomial $y=y(x_1,\dots,x_n)$, let $N(y)$ denote the
\emph{Newton polytope} of~$y$.
That is, $N(y)$ is the convex hull in $\RR^n$ of all lattice points
$\mathbf{m}=(m_1, \dots, m_n)$ such that the coefficient of the monomial
$\mathbf{x}^\mathbf{m}=x_1^{m_1} \cdots x_n^{m_n}$ in $y$ is nonzero.
We claim that the desired valuation $\nu$ can be defined by
\[
\nu (y) =
\min_{\mathbf{m}=(m_1,\dots,m_n)\in N(y)}
(m_1 \nu_1+\cdots+ m_n\mu_n)\,.
\]
The properties (\ref{eq:valuation-0})--(\ref{eq:valuation-2}) are
obvious from this definition,
as is the fact that $\nu(x_i) = \nu_i$ for all~$i$.
The remaining property~(\ref{eq:valuation-2'})
is an immediate consequence of the following statement:
for $x,y \in \overline\AA_{\rm sf}\,$,
the Newton polytope $N(x+y)$ is the convex hull of the union $N(x)\cup N(y)$.
Clearly, $N(x+y)$ is always contained in the convex hull of $N(x)\cup N(y)$.
Furthermore, each vertex $\mathbf{m}$ of this convex hull is either a vertex of
$N(x)$ that lies outside~$N(y)$, or a vertex of $N(y)$ that lies
outside~$N(x)$, or a vertex of both $N(x)$ and~$N(y)$.
We only need to show that in the latter case, the coefficient of
$\mathbf{x}^\mathbf{m}$ in $x+y$ is nonzero; this is clear
since the subtraction-free condition implies that the coefficient of $\mathbf{x}^\mathbf{m}$ in
each of $x$ and $y$ is a nonzero subtraction-free rational expression
in the elements of~$\PP$.
\end{proof}

\noindent
\textbf{Proof of Proposition~\ref{pr:acyclic-necessary-1}.}
It is enough to show that
there exists a tropical valuation $\nu$ on $\AA$ with the
following properties:
\begin{align}
\label{eq:valuation-3}
&\nu(x_i) > 0 \,\, {\rm for} \,\, i=1,\dots, \ell; \\
\label{eq:valuation-4}
&\nu(x_i) = 0 \,\, {\rm for} \,\, i=\ell+1,\dots, n; \\
\label{eq:valuation-5}
&\nu(x'_i) \geq 0 \,\, {\rm for} \,\, i=1,\dots, \ell; \\
\label{eq:valuation-6}
&\nu(y) < 0.
\end{align}
Indeed, in view of (\ref{eq:valuation-0})--(\ref{eq:valuation-2})
and (\ref{eq:valuation-3})--(\ref{eq:valuation-4}),
the exchange relations imply that
$\nu(x'_i) \geq 0$ for $i=\ell+1, \dots, n$,
so $\nu$ takes nonnegative values at all generators
$x_i$ and $x'_i$ of~$\Lower (\Sigma)$.
Then the conditions (\ref{eq:valuation-0})--(\ref{eq:valuation-2})
imply that
$\nu(x) \geq 0$ for all nonzero $x \in \Lower (\Sigma)$.
Comparing this with (\ref{eq:valuation-6}), we conclude that
$y \notin \Lower (\Sigma)$, as claimed.

To prove the existence of a tropical valuation satisfying
(\ref{eq:valuation-3})--(\ref{eq:valuation-6}), we proceed by
induction on $\ell$.
The base is $\ell = 3$.
In this case, the matrix entries $b_{12}, b_{23}$ and $b_{31}$ are
positive while $b_{21}, b_{32}$ and $b_{13}$ are negative.
We set $\nu_1 = \min(|b_{21}|, |b_{31}|)$, $\nu_2=\nu_3=1$, and
$\nu_i=0$ for $3 < i \leq n$.
Let $\nu$ be a tropical valuation with $\nu(x_i) = \nu_i$
for all~$i$, as provided by Lemma~\ref{lem:valuation}.
It obviously satisfies (\ref{eq:valuation-3}) and~(\ref{eq:valuation-4}).
Using the exchange relations (\ref{eq:exchange}) and applying
(\ref{eq:valuation-0}), (\ref{eq:valuation-1}) and (\ref{eq:valuation-2'}),
we also see that
\begin{align*}
\nu(x'_1) &= \min(|b_{21}|\nu_2, |b_{31}|\nu_3)-\nu_1 = 0,\\
\nu(x'_2) &= \min(|b_{12}|\nu_1, |b_{32}|\nu_3)-\nu_2 \geq 0,\\
\nu(x'_3) &= \min(|b_{13}|\nu_1, |b_{23}|\nu_2)-\nu_3 \geq 0.
\end{align*}
Thus, (\ref{eq:valuation-5}) also holds.
To verify (\ref{eq:valuation-6}), note that $y$ is obtained by
exchanging $x_2$ from the cluster $\xx^{(1)} = \{x'_1, x_2, x_3, \dots, x_n\}$.
This cluster has only two elements (namely, $x_2$ and~$x_3$) with positive
value of $\nu$; the rest of its elements have valuation~$0$.
The exchange relation for $y$ is either of the form
$$y x_2 = M {x'_1}^a + N x_3^b$$
or of the form
$$y x_2 = M + N {x'_1}^a x_3^b,$$
where $a$ and $b$ are nonnegative integers, and
$M$ and $N$ are two monomials in
$x_4, \dots, x_n$ with coefficients in~$\PP$.
In both cases, $\nu(M)=\nu(N)=0$, so the first monomial on the right-hand side has
valuation~$0$, while the second one has nonnegative valuation.
It then follows from (\ref{eq:valuation-2'})
that $\nu (y x_2) = 0$, hence $\nu(y) = - \nu(x_2) = -1$, as required.

For the inductive step, assume that $\ell > 3$, and that our
statement holds for the smaller values of $\ell$.
To apply this assumption, we change our choice of an initial seed,
replacing $\Sigma$ with $\Sigma^{(1)}$.
The definition of a matrix mutation implies that
$\Gamma(\Sigma^{(1)})$ has oriented cycle
$2 \to 3 \to \cdots \to \ell \to 2$ (as an induced subgraph).
Applying the inductive assumption to this cycle, we conclude that there
exists a tropical valuation $\nu$ satisfying (\ref{eq:valuation-6})
and the following counterparts of
(\ref{eq:valuation-3})--(\ref{eq:valuation-5}):
\begin{align}
\label{eq:valuation-3'}
&\nu(x_i) > 0 \ \ {\rm for} \,\, i=2,\dots, \ell; \\
\label{eq:valuation-4'}
&\nu(x'_1)=0, \ \ {\rm and}\ \ \nu(x_i) = 0 \ \ {\rm for} \,\,
i=\ell+1,\dots, n; \\
\label{eq:valuation-5'}
&\nu(x''_i) \geq 0 \ \ {\rm for} \,\, i=2,\dots, \ell,
\end{align}
where $x''_i$ is obtained by exchanging $x_i$ from the cluster
$\xx^{(1)} = \{x'_1, x_2, \dots, x_n\}$ of the seed $\Sigma^{(1)}$.
To complete the proof, it suffices to show that this valuation
satisfies (\ref{eq:valuation-3})--(\ref{eq:valuation-5}).
By the mutation rules, we can assume that $x''_i = x'_i$ for $i = 3, \dots, \ell-1$.
Thus, it remains to show that (\ref{eq:valuation-3'})--(\ref{eq:valuation-5'})
imply that
$$\nu(x_1) > 0, \quad \nu(x'_2) \geq 0, \quad \nu(x'_\ell) \geq 0.$$

The exchange relation between the seeds $\Sigma$ and $\Sigma^{(1)}$
can be written as
$$x_1 x'_1 = M x_2^{|b_{21}|} + N x_\ell^{|b_{\ell 1}|},$$
where $M$ and $N$ are monomials in $x_{\ell+1}, \dots, x_n\,$.
Hence $\nu(M) = \nu(N) = 0$, and
$$\nu(x_1) = \nu(x_1 x'_1) = \min(|b_{21}| \nu(x_2), |b_{\ell 1}|
\nu(x_\ell)) > 0,$$
as required.
Similarly,
\begin{equation}
\label{eq:x'_2}
\nu(x_2 x'_2) = \min(|b_{12}| \nu(x_1), |b_{32}| \nu(x_3)) =
\min(|b_{12}b_{21}| \nu(x_2), |b_{\ell 1} b_{12}| \nu(x_\ell), |b_{32}| \nu(x_3)).
\end{equation}
On the other hand, using the mutation rule, we can write down
the exchange relation between $\Sigma^{(1)}$ and $\Sigma^{(2)}$ as follows:
$$x_2 x''_2
= M' {x'_1}^{|b_{12}|} x_3^{|b_{32}|} + N' x_\ell^{|b_{\ell 1}b_{12}|},$$
where, as before, $\nu(M') = \nu(N') = 0$.
Therefore,
\begin{equation}
\label{eq:x''_2}
\nu(x_2 x''_2) =
\min(|b_{\ell 1} b_{12}| \nu(x_\ell), |b_{32}| \nu(x_3)).
\end{equation}
Comparing (\ref{eq:x'_2}) and (\ref{eq:x''_2}), we obtain
$$\nu(x'_2) = \min(\nu(x''_2), (|b_{12}b_{21}|-1) \nu(x_2)) \geq 0,$$
as required.
The remaining inequality $\nu(x'_\ell) \geq 0$ is proved in
the same way, replacing the indices $2$ and $3$ with $\ell$ and
$\ell-1$, respectively.
This completes the proof.
\qed

\section{Proofs of Theorem~\ref{th:fin-gen=tot-cyclic}
and Proposition~\ref{pr:Conway}}
\label{sec:proof-of-skew-symm-tot-cyclic}

\noindent
\textbf{Proof of Theorem~\ref{th:fin-gen=tot-cyclic}.}
The ``if'' part is immediate from Theorem~\ref{th:acyclic-sufficient-lower}.
It remains to prove the converse.

Let $\AA$ be a 
skew-symmetrizable cluster algebra of
rank~3 which has no acyclic seeds.
We need to show that it is \emph{not} finitely generated.
As in \cite{fz-clust1}, we work on the universal cover of the
exchange graph of~$\AA$, i.e., on the $3$-regular tree $\TT_3$ whose
edges are colored in a proper way using the colors $1,2,3$.
As in \cite{fz-clust1}, we write $t \overunder{j}{} t'$ if vertices
$t,t'\in\TT_3$ are joined by an edge labeled by~$j$.
Every vertex $t \in \TT_3$ is associated to a seed
$\Sigma(t)$, and we denote its cluster variables by $x_i (t)$, for $i =1,2,3$.
Since no seed is acyclic, the exchange relation (\ref{eq:exchange})
corresponding to an edge $t \overunder{j}{} t'$ has the form
\begin{equation}
\label{eq:exchange-3}
x_j(t) x_j(t')
= p^+_j x_i(t)^{|b_{ij}|}+ p^-_j x_k(t)^{|b_{kj}|},
\end{equation}
where $(i,j,k)$ is a permutation of~$(1,2,3)$, and both exponents
$|b_{ij}|$ and $|b_{kj}|$ are positive integers.

We fix a root vertex~$t_\circ\in\TT_3\,$.
Let $\nu$ be a \emph{tropical valuation} of $\overline\AA$ in the sense of
Definition~\ref{def:tropical-valuation}.
We abbreviate $\nu_i(t) = \nu (x_i(t))$.
Recall that, by Lemma~\ref{lem:valuation},
the values $\nu_i(t_\circ)$ at the root vertex $t_\circ$ can be assigned
arbitrarily.
The rest of the values $\nu_i(t)$ are then determined by the tropical
version of~(\ref{eq:exchange-3}):
\begin{equation}
\label{eq:exchange-tropical-rk3}
\nu_j (t) + \nu_j (t') =
\min(|b_{ij}(t)| \nu_i(t), |b_{kj}(t)| \nu_k(t)) .
\end{equation}

Let $d(t,t')$ denote the distance between vertices $t$ and $t'$ in~$\TT_3\,$.
The theorem will follow if we show that,
for any positive integer~$r$, there exists a
tropical valuation $\nu$ such that:

\begin{itemize}

\item $\nu_i(t) \geq 0$ for all
$i$ and all $t \in \TT_3$ with $d(t_\circ,t) \leq r$.

\item $\nu_i(t) < 0$ for some $i$ and some $t \in \TT_3$ with $d(t_\circ,t) = r+1$.

\end{itemize}
Indeed, this would imply that no finite set of cluster
variables generate~$\AA$.
As a consequence, $\AA$ does not have a finite set of generators;
otherwise, expressing each generator as a polynomial in
(finitely many) cluster variables, we would obtain a finite set of those
variables as a generating set for~$\AA$.

We now make use of the assumption that $\AA$ is skew-symmetrizable.
Let us denote $s_k (t) = \sqrt{|b_{ij}(t)b_{ji}(t)|}$.
In view of \cite[Lemmas~8.3, 8.4]{fz-clust2},
there are three positive numbers $h_1,h_2,h_3$ such that
\begin{equation}
\label{eq:bij=sh/h}
|b_{ij}(t)| = s_k(t)h_j/h_i
\end{equation}
for all $t$.
We also have
$s_j(t) + s_j(t') = s_i(t) s_k(t)$
for every edge $t \overunder{j}{} t'$ in $\TT_3$.
Renormalizing the tropical evaluations $\nu_j(t)$ by
\begin{equation}
\label{eq:nu-renormalized}
\nu_j(t) = h_j\, s_j(t)\, \delta_j (t),
\end{equation}
and substituting (\ref{eq:bij=sh/h}) and~(\ref{eq:nu-renormalized})
into~(\ref{eq:exchange-tropical-rk3}),
we see that the numbers $\delta_i(t)$ satisfy
\begin{equation}
\label{eq:sjt-delta}
\begin{array}{rcl}
s_j(t)\delta_j(t)+s_j(t')\delta_j(t')
&=&s_k(t) s_i(t) \min(\delta_i(t), \delta_k(t))\\
&=&(s_j(t) + s_j(t'))\min(\delta_i(t), \delta_k(t)).
\end{array}
\end{equation}
Introducing the notation
\[
u_j(t) 
= \frac{s_j(t)}{s_j(t) + s_j(t')}\,,
\]
we rewrite \eqref{eq:sjt-delta} as
\begin{equation}
\label{eq:exchange-tropical-rk3-delta}
u_j(t) \delta_j (t) + u_j(t') \delta_j (t') =
\min(\delta_i(t), \delta_k(t)).
\end{equation}
Note that $u_j(t) + u_j(t') = 1$.
It remains to prove the following lemma.

\begin{lemma}
\label{lem:delta-negative}
Suppose that every edge $t \overunder{j}{} t'$ in $\TT_3$
is assigned two positive numbers $u_j(t)$ and $u_j(t')$
such that $u_j(t) + u_j(t') = 1$.
Then, for any integer $r\geq 0$, there exists a family of real
numbers $\delta_j(t)$ satisfying {\rm (\ref{eq:exchange-tropical-rk3-delta})}
and such that
\begin{itemize}

\item $\delta_i(t) \geq 0$ for all
$i$ and all $t \in \TT_3$ with $d(t_\circ,t) \leq r$.

\item $\delta_i(t) < 0$ for some $i$ and some $t \in \TT_3$ with
 $d(t_\circ,t) = r+1$.

\end{itemize}
\end{lemma}

\begin{proof}
A solution of the relations
(\ref{eq:exchange-tropical-rk3-delta})
is uniquely determined by the initial terms
$\delta_1(t_\circ),\delta_2(t_\circ),\delta_3(t_\circ)$
which can be chosen arbitrarily.
Let us pick these terms so that
$\delta_1 (t_\circ) = \delta_2 (t_\circ) < \delta_3 (t_\circ)$.
It is easy to verify by induction on $d(t,t_\circ)$
that whenever $t \overunder{j}{} t'$ and
$d(t',t_\circ)>d(t,t_\circ)$,
we have $\delta_j (t') \leq \delta_j(t)$; moreover,
this inequality is strict everywhere outside the path going through
$t_\circ$ that has alternating colors
$1$ and~$2$.
It follows that the sequence $\delta(0), \delta(1), \dots$ defined by
\[
\delta(r) = \min(\delta_i(t): i\in\{1,2,3\},\ d(t,t_\circ)=r)
\]
is strictly decreasing.
Now for any given $r$, we obtain the desired solution of
(\ref{eq:exchange-tropical-rk3-delta}) by simply subtracting $\delta(r)$
from all~$\delta_j(t)$.
\end{proof}

\noindent
\textbf{Proof of Proposition~\ref{pr:Conway}.}
Every matrix mutation transforms the matrix
\[
B=\left[\begin{array}{ccc}
0  &  2  & -2\\
-2 & 0 & 2\\
2  & -2  & 0 \\
\end{array}\right]
\]
into $-B$.
Therefore any seed that is mutation equivalent
to $\Sigma$ has exchange matrix~$\pm B$.
It follows that all such seeds are coprime,
hence $\overline\AA = \Upper(\Sigma)$ by
Corollary~\ref{cor:upper-bound-universal}.

Since both exponents $|b_{ij}|$ and $|b_{kj}|$
in each exchange relation (\ref{eq:exchange-3}) are equal to~$2$,
it follows from (\ref{eq:exchange-tropical-rk3}) that there is a
tropical valuation $\nu$ on $\overline\AA$ such that
$\nu(x) = 1$ for each cluster variable~$x$.
This valuation actually makes $\AA(\Sigma)$ into a
graded algebra with the zero degree component~$\ZP$.
Thus, to prove the proposition, it is enough to construct a
non-constant element $y \in \Upper(\Sigma)$ with $\nu(y) = 0$.

The exchange relations from the initial seed $\Sigma$
are of the form
\begin{align*}
&x_1 x'_1 =
 p^-_1 x_2^2+ p^+_1 x_3^2; \\
&x_2 x'_2 =
 p^-_2 x_3^2+ p^+_2 x_1^2; \\
&x_3 x'_3 =
 p^-_3 x_1^2+ p^+_3 x_2^2.
\end{align*}
Then
\[
y \stackrel{\rm def}{=}
\frac{p_1^+ p_2^+ x_1^2 + p_1^- p_2^- x_2^2
+ p_1^+ p_2^- x_3^2}{x_1 x_2}
= \frac{p_1^+ p_2^+ x_1 + p_2^- x'_1}{x_2}
= \frac{p_1^- p_2^- x_2 + p_1^+ x'_2}{x_1}\in \Upper(\Sigma)\,.
\]
Since $y$ is a non-constant element of degree $0$, we are done.
\qed

\section*{Acknowledgments}

Part of this work was done while two of the authors (S.F. and A.Z.) were visiting
the University of Sydney in July 2002.
We thank Gus Lehrer and Joost van Hamel for their kind
hospitality.

\end{document}